\documentclass[11pt]{article}
\usepackage{amsfonts}
\usepackage{amsmath}
\usepackage{amssymb}
\usepackage{amsthm}
\usepackage{array}
\usepackage{bbm}
\usepackage{color}
\usepackage{enumerate}
\usepackage[mathscr]{euscript}
\usepackage{graphicx}
\usepackage[utf8]{inputenc}
\usepackage{latexsym}
\usepackage{mathrsfs}
\usepackage{comment}
\usepackage{mathtools}
\usepackage{appendix}
\usepackage{multirow}
\usepackage{subfigure}
\usepackage[a4paper, top=2.54cm, bottom=2.54cm, left=2.54cm, right=2.54cm]{geometry}
\usepackage{tikz}
\usepackage[font=small]{caption}
\usepackage[square,sort,comma,numbers]{natbib}
\RequirePackage[colorlinks,citecolor=blue,urlcolor=blue,hypertexnames=false]{hyperref}
\usepackage[american]{babel}
\newcommand{\R}{\mathbb{R}}
\newcommand{\E}{\operatorname{\mathbb{E}}}
\newcommand*{\fatdot}[1]{\pmb{\cdot}}
\newcommand{\diff}{\mathrm{d}}

\newcommand{\A}{\mathcal{A}}

\newcommand{\F}{\mathcal{F}}

\newcommand{\be}{\begin{equation}}
\newcommand{\ee}{\end{equation}}
\newcommand{\bea}{\begin{eqnarray}}
\newcommand{\eea}{\end{eqnarray}}
\newcommand{\beas}{\begin{eqnarray*}}
\newcommand{\eeas}{\end{eqnarray*}}
\newcommand{\ds}{\displaystyle}

\def\theequation{\thesection.\arabic{equation}}

\newtheorem{theorem}{Theorem}

\newtheorem{corollary}[theorem]{Corollary}

\newtheorem{definition}[theorem]{Definition}
\newtheorem{example}[theorem]{Example}

\newtheorem{lemma}[theorem]{Lemma}

\newtheorem{proposition}[theorem]{Proposition}
\newtheorem{remark}[theorem]{Remark}

\DeclareMathOperator*{\esssup}{ess\,sup}
\DeclareMathOperator*{\essinf}{ess\,inf}

\begin{document}
\title{\vskip -1.8cm 
Geometric BSDEs
} 
\author{Roger J.~A.~Laeven\footnote{Corresponding author.}  \\
{\footnotesize Dept.~of Quantitative Economics}\\
{\footnotesize University of Amsterdam, CentER}\\
{\footnotesize and EURANDOM, The Netherlands}\\
{\footnotesize \texttt{R.J.A.Laeven@uva.nl}}\\
\and Emanuela Rosazza Gianin \\
{\footnotesize Dept.~of Statistics and Quantitative Methods}\\
{\footnotesize University of Milano-Bicocca, Italy}\\
{\footnotesize \texttt{emanuela.rosazza1@unimib.it}}\\
\and Marco Zullino \\
{\footnotesize Dept.~of Mathematics and Applications}\\
{\footnotesize University of Milano-Bicocca, Italy}\\
{\footnotesize \texttt{M.Zullino@campus.unimib.it}}}

\date{This Version: August 2, 2026}

\maketitle

\vspace{-0.5cm}
\begin{abstract}
We introduce Geometric Backward Stochastic Differential Equations (GBSDEs)
and two-driver BSDEs, which arise naturally in the geometric dynamics
of dynamic return risk measures and of recursive portfolio choice. 
Through
a reduction to auxiliary ordinary BSDEs with logarithmic and singular quadratic (LN-Q) growth rate
$y|\ln (y)|+|z|^2/y$, we establish existence, regularity, uniqueness and
stability of solutions, under both bounded and unbounded driver coefficients and
terminal conditions, and we transfer these results to the original two-driver
equations. We then deploy the theory in two applications. We solve a
portfolio optimization problem under stochastic differential utility,
in which the opportunity process satisfies an endogenously derived
two-driver BSDE and optimality is established via our two-driver
comparison theorem. We further apply GBSDEs to dynamic return and
star-shaped risk measures, including (robust) $L^p$-norms, and
characterize their positive homogeneity, 
star-shapedness and multiplicative convexity.
\\[2mm]
\noindent \textbf{OR/MS Classification:} Probability: Diffusion; Probability: Stochastic model applications; Decision analysis: Risk.
\end{abstract}


\section{Introduction}
The geometric nature of random growth and risky asset price processes in continuous time is well recognized and omnipresent in stochastic analysis and its applications, with the geometric Brownian motion as the canonical example.
By contrast, continuous-time robust dynamic risk assessment via Backward Stochastic Differential Equations (BSDEs) usually occurs in an arithmetic environment (see, e.g., \cite{P97,EPQ97,HIM05,MS05,J08,LS14,KLLSS18,LSSS24}).

In the recent literature (see, e.g., \cite{BLR18,BLR21,LR22,LRZ24} and also \cite{AS08,FH09,DK12,BZ16}), a growing interest has been devoted to return risk measures that assess relative financial positions (or log returns) instead of absolute positions as with monetary risk measures (\cite{FS11,D12}). 
The relative evaluation, reminiscent of the classical notion of relative risk aversion (\cite{P64}),\footnote{That is, monetary and return risk measures exhibit a similar relationship as absolute and relative risk aversion, arithmetic and geometric means, arithmetic and geometric growth, and arithmetic and geometric Brownian motion.} naturally leads to a multiplicative, i.e., geometric, structure for return risk measures and their acceptance sets.
These considerations --- together with the fact that a wide family of dynamic monetary risk measures is induced by BSDEs (see \cite{CHMP,P05,RG06,BN08,DPR10}) --- motivate the study of geometric forms of BSDEs, both from a mathematical point of view as well as for applications such as portfolio optimization. 

The recent \cite{LRZ23} characterizes properties of dynamic return and star-shaped risk measures induced via ordinary one-driver BSDEs under classical Lipschitz-type $L^{2}$- or $L^{\infty}$-standard assumptions, representing star-shaped functionals as minima of convex functionals induced by Lipschitz BSDEs.
However, \cite{LRZ23} does not address the inherently geometric nature of dynamic return risk measures and geometric dynamic risk assessment more generally. 
Its results are confined to the classical, additive Lipschitz 
regime and are fundamentally unrelated to the non-Lipschitz, logarithmic and singular quadratic growth rates that, as we will see below, arise naturally from the geometric structure.
To address these problems, we develop a structurally novel approach by introducing and analyzing the stochastic differential equations that naturally govern the dynamics of return risk measures and geometric risk assessment, establishing their solutions' existence, regularity, uniqueness, and stability, and applying the resulting theory to portfolio optimization and dynamic risk measurement. 

More specifically, we introduce a novel class of BSDEs, which we term \textit{Geometric BSDEs} (GBSDEs). 
Just like the geometric mean is connected to the arithmetic mean via $\mathbb{G}[X]:=\exp(\mathbb{E}[\ln (X)])$, GBSDEs arise naturally from the substitution $\tilde \rho_{t}(X):=\exp(\rho_t(\ln(X)))$ for any $t\in[0,T]$, where $T>0$ is a fixed finite time horizon and $Y_{\fatdot{}}:=\rho_{\fatdot{}}(\ln(X))$ is the first component of the solution $(Y_{\fatdot{}},Z_{\fatdot{}})$ to the quadratic BSDE \
$Y_t=\ln(X)+\int_t^T f(s,Y_s,Z_s)ds-\int_t^TZ_sdW_s.$
Here, $\left(W_{t}\right)_{t\in[0,T]}$ is a standard $n$-dimensional Brownian motion, the driver $f:\Omega\times[0,T]\times\R\times\R^n\to\R$ is $d\mathbb{P}\times dt$-a.s.\ continuous in $(y,z)$ and satisfies $|f(t,y,z)|\leq \alpha_t+\beta_t|y|+\gamma_t|z|+\delta |z|^2$ for suitable positive stochastic coefficients $\alpha,\beta,\gamma$ and a constant $\delta>0$, whereas $X$ is a strictly positive random variable verifying further integrability conditions. 
Under appropriate hypotheses, It\^o's formula yields the dynamics of $\tilde Y_{\fatdot{}}:=\tilde\rho_{\fatdot{}}(X)$ as the first component of the solution to a GBSDE of the form
\begin{equation}
\tilde Y_t=X+\int_t^T\tilde Y_s\tilde f(s,\tilde Y_s,\tilde Z_s)ds-\int_t^T\tilde Y_s\tilde Z_sdW_s,
\label{eq:GBSDEintro}
\end{equation}
where the driver $\tilde{f}(t,y,z):=f(t,\ln(y),z)-\frac{1}{2}|z|^2$ verifies the growth rate
\begin{equation}
|\tilde f(t,y,z)|\leq \alpha_t +\beta_t|\ln(y)|+\gamma_t|z|+\Big(\delta+\frac{1}{2}\Big) |z|^2.
\label{eq:GRintro}
\end{equation}
The multiplicative nature of \eqref{eq:GBSDEintro} is even more evident in its infinitesimal version \[-\frac{d\tilde Y_t}{\tilde Y_t}=\tilde f(t,\tilde Y_t,\tilde Z_t)dt-\tilde Z_tdW_t, \ \ \ \tilde Y_T=X.\]
GBSDEs generalize the \textit{geometric} martingale representation theorem (\cite{BEK05,DR09}), which corresponds to $\tilde f\equiv 0$, and are naturally connected to dynamic return risk measures (\cite{BLR18,BLR21,LR22}).

For a systematic analysis of its solution's existence, regularity, uniqueness and stability, we embed the GBSDE~\eqref{eq:GBSDEintro} into the substantially broader class of \textit{two-driver BSDEs}, which is also used essentially later in our portfolio choice application:
\begin{equation}
Y_t=X+\int_t^Tg_1(s,Y_s,Z_s)ds-\int_t^Tg_2(s,Y_s,Z_s)dW_s,
\label{eq:2DBSDEintro}
\end{equation}
with $g_1:\Omega\times[0,T]\times\R_+\times\R^n\to\R_+$ and $g_2:\Omega\times[0,T]\times\R_+\times\R^n\to\R^n$.
Two-driver BSDEs were introduced in \cite{PP90} under restrictive hypotheses: a bi-Lipschitz condition on $g_2$ in $(y,z)$ with injectivity in $z$, and Lipschitz continuity of $g_1$.
While significant progress has been made to weaken the Lipschitz hypothesis on $g_1$ (see, e.g., \cite{BH06,BEH15,BEO17,DHK13,K00}), little attention has been given to BSDEs with general double drivers, since $g_2(t,y,z)=z$ has been convenient in many applications (\cite{BN08,BEK05}). 
However, the geometric structure in~\eqref{eq:GBSDEintro} yields $g_2(t,y,z)=yz$, for which the bi-Lipschitz assumption fails. 
We relax this condition by requiring only that $|g_2|$ grows sufficiently fast in $(y,z)$, and our general assumptions, including those on $g_1$, ensure that~\eqref{eq:GBSDEintro} is a special case of~\eqref{eq:2DBSDEintro}.

Our results for two-driver BSDEs are in substantial part --- but importantly not fully --- obtained by deriving corresponding results for an auxiliary class of one-driver BSDEs.
Setting $\bar{Z}_{t}:=\tilde{\rho}_{t}\tilde{Z}_{t}$ in~\eqref{eq:GBSDEintro} yields an ordinary BSDE whose driver satisfies the \textit{logarithmic-quadratic} (LN-Q) growth rate
\begin{equation}
       g(t,y,z)\leq \alpha_t y+\beta_t y|\ln(y)|+\gamma_t|z|+\delta|z|^2/y.
\label{eq:lnQgrowthintro}
\end{equation}
Hence, as we will see, an auxiliary (but not sufficient) tool to study~\eqref{eq:2DBSDEintro} is the ordinary BSDE
\begin{equation}
Y_t=X+\int_t^{T}g(s,Y_s,Z_s)ds-\int_t^TZ_sdW_s,
\label{eq:BSDEintro}
\end{equation}
where $g$ satisfies~\eqref{eq:lnQgrowthintro}, exhibiting a logarithmic nonlinearity and quadratic singularity at zero.
The introduction and comprehensive analysis of general two-driver BSDEs, encompassing GBSDEs, and their connections to portfolio choice and dynamic risk measures, constitutes our main methodological contribution.

\medskip
\noindent\textbf{Main results.}
We establish existence, regularity, uniqueness and stability results for~\eqref{eq:2DBSDEintro}, and along the way \eqref{eq:lnQgrowthintro}--\eqref{eq:BSDEintro}, allowing both bounded and unbounded stochastic driver growth coefficients $\alpha,\beta,\gamma$ 
and terminal conditions $X$.
While existence, uniqueness and stability for BSDEs with $y|\ln(y)|$-growth were established in \cite{BKK17}, and existence and uniqueness under $|z|^2/y$-growth in \cite{BT19} (without stability), this is the first work to simultaneously treat both nonlinearities in~\eqref{eq:lnQgrowthintro} and to address two-driver BSDEs of the form~\eqref{eq:2DBSDEintro}.
Our assumptions differ from \cite{BKK17} and relax those in \cite{BT19}.
They are fundamentally unrelated to the Lipschitz standard assumptions considered in \cite{LRZ23}, see also \cite[Chapter 4]{Z17}, and intrinsically different from classical bi-Lipschitz two-driver BSDEs (\cite{PP90}) and standard quadratic (F)BSDEs (\cite{K00,BH06}, \cite[Chapter 7]{Z17}, \cite{LMO22}). 
Indeed, unlike standard quadratic (F)BSDEs, which retain the martingale term \(Z\,dW\), our framework includes a nonlinear martingale driver \(g_2\) that is connected (via \eqref{eq:lnQgrowthintro}) to singular LN-Q growth
\(y|\ln (y)|+|z|^2/y\), adding genuinely new features compared to the existing quadratic framework.

The existence proof draws some inspiration from the insightful superlinear analysis in \cite{B20,BEH15}. 
The regularity of the $z$-component requires novel \textit{a priori} estimates: square-integrability under $|z|^2/y$-growth was previously known only for bounded terminal conditions bounded away from zero (\cite{BT19}), and we remove both restrictions.
For uniqueness, we assume joint convexity of $g$ in $(y,z)$ and employ a $\theta$-difference strategy in the spirit of \cite{BH08}, but modified to the $y|\ln(y)|+|z|^2/y$-growth via a suitable generalization of a stochastic Bihari--Gr\"onwall inequality (\cite{HKK22}). 
Unlike \cite{BH08,DHR11,FHT23}, we require neither exponential moments of all orders nor monotonicity in $y$.
Our stability results rely on the comparison principle we establish, and provide as a special case the first stability result under singular quadratic growth.
The transfer to two-driver BSDEs requires a separate analysis for the $z$-regularity, which \emph{cannot} be derived from the single-driver equation.
Comparison and stability results for two-driver BSDEs are also established.

\medskip
\noindent\textbf{Applications.}
We show that the two-driver BSDE framework arises naturally in portfolio optimization under stochastic differential utility (SDU) preferences \`a la \cite{DE92}, where the multiplicative structure of the SDU's opportunity process leads to a genuine two-driver BSDE.
The well-posedness theory developed in this paper for two-driver BSDEs is essential to verify optimality of the candidate portfolio strategy, an application that inherently requires the full framework and cannot be handled by reducing it to a one-driver BSDE with growth rate already known in the literature.

Furthermore, we apply the theory to dynamic return and star-shaped risk measures induced by GBSDEs, defining them under natural and weak conditions, 
and identifying sufficient conditions on the drivers for positive homogeneity, star-shapedness and multiplicative convexity.
In particular, we show that $L^{p}$-norms and robust $L^{p}$-norms --- the return counterparts of entropic and robust entropic risk measures (\cite{FS11,LS13}) --- can be described as solutions to GBSDEs, which to our knowledge is the first such representation.

\medskip
\noindent\textbf{Outline.}
Section~\ref{sec:prel} provides the basic framework and preliminaries. 
Section~\ref{sec:gbsdeLNQ} introduces GBSDEs and presents results for the LN-Q BSDE~\eqref{eq:lnQgrowthintro}--\eqref{eq:BSDEintro}. 
Section~\ref{sec:2DBSDEs} develops the theory of two-driver BSDEs.
Section~\ref{sec:app_lnq_control} develops the portfolio application with a log-type SDU, whereas
Section~\ref{sec:appRRM} applies GBSDEs to dynamic return and star-shaped risk measures.
Auxiliary results and proofs not contained in the main text are provided in Online Appendix~\ref{app}.
\setcounter{equation}{0}

\section{Preliminaries}\label{sec:prel}
In this section, we introduce the notation and setting utilized in the sequel. 
Let $\left( \Omega ,\mathcal{F},\mathbb{P}\right)$ be a probability space. 
Furthermore, let $T>0$ denote a finite time horizon and let $\left( W_{t}\right) _{t\in[0,T]}$ be a standard $n$-dimensional Brownian motion defined on $\left( \Omega ,\mathcal{F},\mathbb{P}\right)$. 
We equip the probability space with $\left( \mathcal{F}_{t}\right) _{t\in[0,T]}$, the augmented filtration associated to that generated by $\left(W_{t}\right)_{t\in[0,T]}$. 
We require w.l.o.g.\ that $\mathcal{F}=\mathcal{F}_T$. 
All equalities and inequalities between random variables are understood to hold $d\mathbb{P}$-almost surely whereas for stochastic processes they are meant to be valid $d\mathbb{P}\times dt$-a.s., unless specified otherwise.
We equip the space $L^0$ of all $\F$-measurable random variables with the usual pointwise partial order relation, writing $X\geq Y$ when $X(\omega)\geq Y(\omega)$ for $\mathbb{P}$-almost every $\omega\in\Omega.$  
We use the notation $L^0_{+}$ to denote the set of all strictly positive random variables. 
Given $\mathcal{X}\subseteq L^0$, we set $\mathcal{X}_{+}:=\mathcal{X}\cap L^0_{+}$.
For each fixed $t\in[0,T]$, given $\mathcal{X}\subseteq L^0$, we denote by $\mathcal{X}(\mathcal{F}_t)$ the space of $\mathcal{F}_t$-measurable random variables belonging to $\mathcal{X}.$
For any $p\in[1,+\infty)$, $L^p(\mathcal{F}_T)$ is the set of $p$-integrable random variables whose norm is denoted by $\|\fatdot{}\|_p$. 
$L^{\infty}(\mathcal{F}_t)$ is the space of $\mathcal{F}_t$-measurable and essentially bounded random variables, whose norm is denoted by $\|\fatdot{}\|_{\infty}$.
When no confusion can arise, we write $L^p$ instead of $L^{p}(\mathcal{F}_T)$, and similarly for the other spaces, without further specification. 
For any $x,y\in\R^n$ with $n\in\mathbb{N}$, $n>1$, we write $x\cdot y$ for the usual scalar product in $\R^n$, i.e., $x\cdot y:=\sum_{i=1}^{n}x_iy_i$. 
For clarity, we sometimes use the notation $x\cdot y$ even when $x$ and $y$ are scalars. 
We let $\R_+:=\{x\in\R:x>0\}$.

Next, we define the primary functional spaces that we consider. 
From now on, we will use ``p.p.'' to denote any predictable process with respect to \mbox{$(\mathcal{F}_t)_{t\in[0,T]}\otimes\mathcal{B}(0,T)$}, with $\mathcal{B}$ the Borel $\sigma$-algebra, 
and valued in $\R^n$, with $n\geq1$. 
We let
\begin{eqnarray*}
\hspace{-0.6cm}&&\mathcal{H}^{p}_T := \left\{(Y_t)_{t\in[0,T]} \text{ p.p.: } \mathbb{E}\left[\esssup_{t\in[0,T]}|Y_t|^p\right]<+\infty\right\},\\
	\hspace{-0.6cm}&&\mathcal{H}^{\infty}_T := \left\{(Y_t)_{t\in[0,T]} \text{ p.p.: } \esssup_{(t,\omega)\in[0,T]\times\Omega}|Y_t|<+\infty\right\},\\
 \hspace{-0.6cm}&&\mathcal{L}^{2}_T := \left\{(Y_t)_{t\in[0,T]} \text{ p.p.: } \int_0^T|Y_t|^2dt<+\infty, \ \mathbb{P}\mbox{-a.s.}\right\}, \\
\hspace{-0.6cm}&&\mathcal{M}^{p}_T := \left\{(Y_t)_{t\in[0,T]} \text{ p.p.: } \mathbb{E}\left[\left(\int_0^T|Y_t|^2dt\right)^{\frac{p}{2}}\right]<+\infty\right\},\\
	\hspace{-0.6cm}&&\mbox{BMO}({\mathbb{P}}) := \left\{(Y_t)_{t\in[0,T]} \text{ p.p.: } \exists C>0 \text{ s.t. } \mathbb{E}\left[\left. \int_t^T|Y_s|^2ds \right|\mathcal{F}_t\right]\leq C \ d\mathbb{P}\times dt\text{-a.s.}\right\}. 
\end{eqnarray*}
Sometimes we will employ a different probability measure $\mathbb{Q}$ on $(\Omega,\mathcal{F}).$ 
In that case, we will specify the regularity of a process w.r.t.\ this measure by writing, e.g., $\mathcal{H}^p_T(\mathbb{Q})$ and analogously for the other spaces. 
Similarly, we will specify that the expectations are taken w.r.t.\ $\mathbb{Q}$ by writing $\mathbb{E}_{\mathbb{Q}}\left[\fatdot{}\right].$
Henceforth, we use the shorthand notation $X_{\fatdot{}}$ instead of $(X_t)_{t\in[0,T]}$ for a stochastic process. 
For any $X_{\fatdot{}}\in\mathcal{H}^{\infty}_T$, we define the norm $\|X\|^T_{\infty}:=\ds\esssup_{(t,\omega)\in[0,T]\times\Omega}|X_t|$.  $\mathcal{H}^{\infty}_T$ is a Banach space if equipped with $\|\fatdot{}\|_{\infty}^T$. 
When dealing with $X_{\fatdot{}}\in\mbox{BMO}(\mathbb{P})$, we consider the norm $\|X\|_{\mbox{BMO}}:=\sup_{t\in[0,T]}\mathbb{E}[\int_t^T|X_s|^2ds|\mathcal{F}_t]$. 
As is well known, \mbox{$(\mbox{BMO}(\mathbb{P}),\|\fatdot{}\|_{\mbox{BMO}})$} forms a Banach space (\cite{K94}).

We define the stochastic exponential of  $\gamma\in\mathcal{L}_T^2$ as 
$\mathcal{E}^{\gamma}_t:=\exp\left(\int_0^t\gamma_sdW_s-\frac{1}{2}\int_0^t|\gamma_s|^2ds\right).$
If $\gamma\in\text{BMO}(\mathbb{P})$, $\mathcal{E}_T^{\gamma}$ is the density of an equivalent probability measure $\mathbb{Q}^{\gamma}$ (w.r.t.\ $\mathbb{P}$) and $W^{\gamma}_{\fatdot{}}:=W_{\fatdot{}}-\int_0^{\fatdot{}}\gamma_sds$ is a $\mathbb{Q}^{\gamma}$-Brownian motion w.r.t.\ the same stochastic base of $W_{\fatdot{}}$ (by Girsanov).

In what follows, we also need the space of positive random variables bounded away from zero. 
Specifically, fixing $\varepsilon>0$, we introduce the set and respective functional space 
$$\mathcal{L}^{\infty}_{\varepsilon}(\mathcal{F}_T):=\{Y\in L^{\infty}(\mathcal{F}_T): Y\geq\varepsilon \mbox{ a.s.}\},
\qquad\mathcal{L}^{\infty}(\mathcal{F}_T):=\ds\bigcup_{\varepsilon>0}\mathcal{L}^{\infty}_{\varepsilon}(\mathcal{F}_T).$$ 
When no confusion is possible, we simply write $\mathcal{L}^{\infty}_{\varepsilon}$ and $\mathcal{L}^{\infty}$, respectively. 
In addition, when working with a filtered probability space $(\Omega,\mathcal{F},(\mathcal{F}_t)_{t\in[0,T]},\mathbb{P})$, we write $\mathcal{L}^{\infty}(\mathcal{F}_t)$ (resp.\ $\mathcal{L}^{\infty}_{\varepsilon}(\mathcal{F}_t)$) to indicate the set of random variables in $\mathcal{L}^{\infty}$ (resp.\ $\mathcal{L}^{\infty}_{\varepsilon}$) that are $\mathcal{F}_t$-measurable. 
See Lemma~\ref{lem:spaceidentification} for a characterization of $\mathcal{L}^{\infty}$, which will often be used in the following. 
We will also make use of the following two related spaces for stochastic processes: 
\begin{align*}
        \mathcal{L}^{\infty}_{\varepsilon,T}:= \left\{Y\in\mathcal{H}^{\infty}_T: \ Y_t\geq\varepsilon \ d\mathbb{P}\times dt\text{-a.s.}\right\}, 
        \qquad\mathcal{L}^{\infty}_{T}:= \bigcup_{\varepsilon>0}\mathcal{L}^{\infty}_{\varepsilon,T}.
\end{align*}

Hereafter, we will consider BSDEs of the form:
\begin{equation}
    Y_t=X+\int_t^Tf(s,Y_s,Z_s)ds-\int_t^TZ_sdW_s.
    \label{eq:BSDEq1}
\end{equation}
The pair $(X, f)$ will be referred to as the \emph{parameters} of the associated BSDE.
\begin{definition}
The couple $(Y,Z)$ is a solution to Equation~\eqref{eq:BSDEq1} if it verifies this equation in the It\^o sense, $Y$ is predictable with continuous trajectories, $Z\in\mathcal{L}_T^2$ is predictable, and $\int_0^T|f(s,Y_{s},Z_{s})|ds<+\infty,$ $d\mathbb{P}$-a.s.  
\label{def:solBSDE}
\end{definition}

\setcounter{equation}{0}
\section{Main Results on LN-Q BSDEs}\label{sec:gbsdeLNQ}
\subsection{From BSDEs to GBSDEs and LN-Q BSDEs}\label{sec:1-1}
The dynamics of return risk measures arise naturally by exponentiating the solution of a BSDE associated with an ordinary dynamic risk measure; see also Appendix~\ref{sec:drdm} for further details. 
More precisely, let \(X>0\), and let \(\rho_t(\ln X)\) denote the first component of the solution to a quadratic BSDE with terminal condition \(\ln X\) and driver \(f\) satisfying
\(
|f(t,y,z)|\le C\bigl(1+|y|+|z|^2\bigr)\) for some \(C>0\).
Defining \(\tilde\rho_t(X):=\exp\bigl(\rho_t(\ln X)\bigr)\), applying It\^o's formula and using
\(
d\rho_t=-f(t,\rho_t,Z_t)\,dt+Z_t\,dW_t
\)
yields
\[
d\tilde\rho_t
=\tilde\rho_t\,d\rho_t+\frac12\tilde\rho_t|Z_t|^2\,dt
=\tilde\rho_t\left(-f(t,\rho_t,Z_t)+\frac12|Z_t|^2\right)dt
+\tilde\rho_t Z_t\,dW_t.
\]
Since \(\rho_t=\ln(\tilde\rho_t)\), setting
\begin{equation}
\tilde f(t,y,z):=f(t,\ln y,z)-\frac12|z|^2\label{eq:fdef}
\end{equation}
shows that \((\tilde\rho,\tilde Z)\) satisfies the multiplicative BSDE
\[
-\frac{d\tilde\rho_t}{\tilde\rho_t}
=
\tilde f(t,\tilde\rho_t,\tilde Z_t)\,dt-\tilde Z_t\,dW_t,
\qquad
\tilde\rho_T=X,
\]
where \(\tilde Z_t:=Z_t\). Moreover, the growth condition on \(f\) yields
\[
|\tilde f(t,y,z)|
\le C\bigl(1+|\ln y|+|z|^2\bigr),
\]
thus exhibiting the logarithmic dependence on \(y\) and quadratic dependence on \(z\) characteristic of the geometric setting.

Motivated by the previous observations, we define a GBSDE as the ensuing Backward Stochastic Differential Equation
\begin{equation}
\label{eq:GBSDE} 
\tilde\rho_t=X+\int_t^T\tilde\rho_s \tilde f(s,\tilde\rho_s,\tilde Z_s)ds-\int_t^T\tilde\rho_s\tilde Z_sdW_s,
\end{equation}
with $\tilde f$ given by \eqref{eq:fdef}.
Setting $\bar Z_t:=\tilde\rho_t \tilde Z_t$ transforms \eqref{eq:GBSDE} into the ordinary BSDE
\[\bar\rho_t=X+\int_t^T\bar\rho_s\tilde f(s,\bar\rho_s,\bar Z_s/\bar\rho_s)ds-\int_t^T\bar Z_sdW_s,\]
with the new driver $\bar f(t,y,z):=\tilde f(t,y,z/y)$ satisfying a \textit{logarithmic-quadratic} (LN-Q) growth rate:
$y\bar f(t,y,z)\leq C(y+y|\ln(y)|+|z|^2/y), \ \forall (y,z)\in\R_+\times\R^n.$
Since we want to allow for stochastic and unbounded coefficients, we consider in the sequel the generalized growth rate
\begin{equation*} 
y\bar f(t,y,z)\leq \alpha_ty+\beta_ty|\ln(y)|+\gamma_t|z|+\delta|z|^2/y, \ \forall (y,z)\in\R_+\times\R^n,
\end{equation*}
where $\alpha,\beta,\gamma$ are predictable and positive stochastic processes and $\delta\geq 0$.

Hence, we study the auxiliary ordinary BSDE
\begin{equation}
Y_t=X+\int_t^Tg(s,Y_s,Z_s)ds-\int_t^TZ_sdW_s,  
\label{eq:LNQ}
\end{equation}
where $g:[0,T]\times\Omega\times\R_+\times\R^n\to\R$ is $\mathcal{P}\times\mathcal{B}(\R_+)\times\mathcal{B}(\R^n)/\mathcal{B}(\R)$-measurable with LN-Q growth rate
\begin{equation}
|g(t,y,z)|\leq \alpha_ty+\beta_ty|\ln(y)|+\gamma_t|z|+\delta|z|^2/y, \ d\mathbb{P}\times dt\text{-a.s.} \ \forall(y,z)\in\R_+\times\R^n.
\label{eq:gLNQ}
\end{equation}

The following preliminary result formalizes a one-to-one correspondence between LN-Q and quadratic BSDEs (\cite{K00}, \cite[Chapter 7]{Z17}) under stringent conditions.
Recall the quadratic BSDE
\begin{equation}
Y'_t= X'+\int_t^Tg'(s, Y'_s, Z'_s)ds-\int_t^TZ'_tdW_s,
\label{eq:QG}
\end{equation}
with 
\begin{equation}|g'(t,y,z)|\leq\alpha_t+\beta_t|y|+\gamma_t|z|+\eta|z|^2, \ d\mathbb{P}\times dt\text{-a.s.} \ \forall(y,z)\in\R\times\R^n,
\label{eq:gGQ}
\end{equation}
where $\eta\geq 0$.  
\begin{proposition}
\label{prop:equivalence}
Let $g$ satisfy condition~\eqref{eq:gLNQ}, 
where $\alpha ,\beta ,\gamma $ are predictable and positive stochastic processes.

If Equation~\eqref{eq:LNQ} admits a positive solution with parameters $(X,g)$ where $X\in L^{0}_{+}(\mathcal{F}_T)$, then Equation~\eqref{eq:QG} admits a solution with parameters $(X',g')$, where $g'$ verifies~\eqref{eq:gGQ} with $\eta=\delta+1/2$ and $ X'\in L^{0}(\mathcal{F}_T)$.

Conversely, if Equation~\eqref{eq:QG} admits a solution with parameters $( X', g')$, where $g'$ verifies~\eqref{eq:gGQ} and $ X'\in L^{0}(\mathcal{F}_T)$, then Equation~\eqref{eq:LNQ} admits a positive solution with parameters $(X,g)$, where $g$ satisfies~\eqref{eq:gLNQ} with $\delta=\eta+1/2$ and $X\in L^{0}_{+}$.
\end{proposition}
\begin{corollary}
With the same notation as in Proposition~\ref{prop:equivalence}, let $X\in\mathcal{L}^{\infty}(\mathcal{F}_T)$ and $\alpha,\beta,\gamma\in \mathcal{H}^{\infty}_T$. 
Then, there exist a maximal and minimal solution to Equation~\eqref{eq:LNQ} with parameters $(X,g)$. 
Each solution between the minimal and the maximal verifies $(Y,Z)\in \mathcal{H}^{\infty}_T\times\text{BMO}(\mathbb{P}).$
\label{cor:gen1to1}
\end{corollary}

As is evident from (the proof of) Corollary~\ref{cor:gen1to1}, its restrictive assumptions, 
which will not be satisfied in any of our applications, are pivotal so far. 
Establishing the existence and properties of the BSDE solution becomes intrinsically different when these assumptions are removed, e.g., when the driver growth coefficients are stochastic and unbounded and/or when the terminal condition is unbounded and possibly not even integrable.  
To investigate the existence, regularity, uniqueness and stability of the solution outside the highly restrictive setting of Corollary~\ref{cor:gen1to1}, 
where the introductory one-to-one correspondence 
breaks down and can no longer be exploited, an essential departure from the standard quadratic and classical bi-Lipschitz two-driver settings of \cite{K00} and \cite{PP90} is required:   
we must explore genuinely distinct approaches that can jointly handle the logarithmic nonlinear growth in the state variable and the quadratic singularity at zero. 
These problems are studied in what follows.

\subsection{Assumptions}\label{sec:assumptions}
For the sake of generality and completeness, we will consider subsets of the following assumptions on the driver $g$ and terminal condition $X$ in Equation~\eqref{eq:LNQ}:
\begin{itemize}
    \item[(H1)] $g:\Omega\times[0,T]\times\R_+\times\R^{n}\to\R_+$ is a predictable stochastic process, continuous in $(y,z)$ $d\mathbb{P}\times dt$-a.s., verifying
    \begin{equation*}
       0\leq g(t,y,z)\leq \alpha_t y+\beta_t y|\ln(y)|+\delta|z|^2/y=:h(t,y,z),
    \end{equation*}
    where $\alpha,\beta$ are positive and predictable stochastic processes and $\delta>0$;
    \item[(H1)'] With the same notation as in (H1), the process $g$ now verifies
    \begin{equation*}
       0\leq g(t,y,z)\leq \alpha_t y+\beta_t y|\ln(y)|+\gamma_t|z|+\delta|z|^2/y=:h(t,y,z),
    \end{equation*}
    where $\gamma:[0,T]\times\Omega\to\R$ is a positive and predictable stochastic process;
\item[(H1)''] With the same notation as in (H1), the process $g$ now verifies
    \begin{equation*}
        0\leq g(t,y,z)\leq \alpha_t y+\beta_t y|\ln(y)|+\eta_t\cdot z+\delta|z|^2/y=:h(t,y,z),
    \end{equation*}
    where $\eta:[0,T]\times\Omega\to\R^n$ is a predictable stochastic process;
    \item[(H2)] Let $X$ be a strictly positive $\mathcal{F}_T$-measurable random variable such that
    \begin{equation*}
       \mathbb{E}\left[(1+X^{2\delta+1})^{e^{B}}\exp\left(e^{B}(A+B)\right)\right]<+\infty,
    \end{equation*}
       where $A:=\int_0^T\alpha_tdt$ and $B:=\int_0^T\beta_tdt$. 
    
    \item[(H2)'] With the same notation as in (H2), let $p>1$ such that
    \begin{equation*}
       \mathbb{E}\left[(1+X^{2\delta+1})^{p(e^{B}+1)}\exp\left(p(e^{B}+1)(A+B)\right)\right]<+\infty.
    \end{equation*}
    In addition, we assume $\alpha,\beta\in \mathcal{H}^{q}_T$, with $\frac{1}{q}+\frac{1}{p}=1$.
    \item[(H2)''] With the same notation as in (H2), let $p>1$ and $\gamma$ be as in (H1)' such that
    \begin{equation*}
       \mathbb{E}\left[(1+X^{2\delta+1})^{p(e^{B}+1)}\exp\left(p(e^{B}+1)\left((A+B)+\frac{1}{4\eta}\int_0^T\gamma_t^2dt\right)\right)\right]<+\infty.
    \end{equation*}
\noindent Here, $\eta\!\in\!(\!0,\!\frac{1\wedge (\!p-1\!)}{2}).$ 
Furthermore, there exists $q'\!>\!0$ such that $\mathbb{E}[\int_0^T\!e^{q'\gamma_t}dt]\!<+\infty$ and $\alpha,\beta\in \mathcal{H}^{q}_T$, with $\frac{1}{q}+\frac{1}{p}=1$.
\end{itemize}
We employ one and the same notation $h$ across (H1), (H1)', and (H1)'', with slight abuse of notation; the respective form of $h$ will be specified on a case-by-case basis. 

\subsection{Existence and regularity}
\begin{proposition}
\label{prop:EUgen}
(i) Assume (H1) and (H2). 
Then, Equation~\eqref{eq:LNQ} with driver $h$ as defined in (H1) admits a positive solution such that
$$\ds\sup_{t\in[0,T]}\mathbb{E}\left[Y_t^{(2\delta+1)e^B}\right]<+\infty, \ \text{ and } Z\in\mathcal{L}^2_{T}.$$

(ii) Assume (H1) and (H2)'. 
Then, there exists a positive solution to Equation~\eqref{eq:LNQ} with driver $h$ as defined in (H1), verifying the further regularity $$\mathbb{E}\bigg[\sup_{t\in[0,T]}Y_t^{p(2\delta+1)(e^B+1)}\bigg]<+\infty, \ \text{ and } Z\in\mathcal{M}_T^2.$$

(iii) Assume (H1)' and (H2)''. 
Then, there exists a positive solution to Equation~\eqref{eq:LNQ} with driver $h$ as defined in (H1)', verifying the same regularity as in the case of (H1)+(H2)'.
\end{proposition}

\begin{theorem}
Assume (H1) and (H2), with driver $g$ as in (H1). 
Then, there exists a positive solution to Equation~\eqref{eq:LNQ} in the sense of Definition~\ref{def:solBSDE} verifying $0<Y\leq Y^h$, where $Y^h$ is a solution to Equation~\eqref{eq:LNQ} with driver $h$. 
We have the regularity $0< \ds\sup_{t\in[0,T]}\mathbb{E}\left[Y_t^{(2\delta+1)e^B}\right]<+\infty$ and $Z\in\mathcal{L}^2_{T}$. 
Furthermore, if we assume (H1) and (H2)' (or, alternatively, (H1)' and (H2)''), then the solution has the further regularity $0< \mathbb{E}\bigg[\ds\sup_{t\in[0,T]}Y_t^{p(2\delta+1)(e^B+1)}\bigg]<+\infty$ and $Z\in\mathcal{M}^2_T$. 
\label{th:Egen}
\end{theorem}

\begin{proof}
We first apply Lemma~\ref{lem:bahlali} in Appendix~\ref{app}. 
We choose $X_1=X_2=X$, $g_1\equiv 0$ and $g_2(t,y,z):=h(t,y,z)=\alpha_t y+\beta_t y|\ln(y)|+\delta |z|^2/y$. 
Clearly, $g_1\equiv0\leq g\leq g_2$. 
The BSDE with parameters $(X,0)$ admits a unique positive solution $Y^0$ by Proposition~1.1~(i) in \cite{B20}, while Proposition~\ref{prop:EUgen} ensures the existence of a positive solution $Y^{h}$ to the BSDE with parameters $(X,h)$. 
Moreover, $0<Y^0\leq Y^{h}$, since after a suitable localization $(\tau_n)_{n\in\mathbb{N}}$,
$$\mathbb{E}[Y^{h}_t|\mathcal{F}_t]=\mathbb{E}\left[ \left. X+\int_t^{\tau_n}g_2(s,Y_s^{h},Z^{h}_s)ds \right|\mathcal{F}_t\right]\geq \mathbb{E}[X|\mathcal{F}_t]=Y^0_t,$$
by positivity of $h$. 
Finally, for any $(t,\omega)\in[0,T]\times\Omega$ and $y\!\in[Y^0_t(\omega),Y^{h}_t(\omega)]$,
$$g(t,y,z)\!\leq\! (\alpha_t+\beta_t)+(\alpha_t+\beta_t)(1\!+\!|y|^2)\!+\!\frac{\delta|z|^2}{y}\leq \!(\alpha_t+\beta_t)\!+\!(\!\alpha_t\!+\!\beta_t)(1\!+\!|Y^{g_2}_t(\omega)|^2)\!+\!\frac{\delta|z|^2}{Y^0_t(\omega)},$$
so $g$ satisfies all conditions of Lemma~\ref{lem:bahlali}. 
This yields a solution $(Y,Z)$ with $0 < Y^0 \leq Y \leq Y^{h}$, hence the desired regularities for $Y$. 
The cases (H1)+(H2)' and (H1)'+(H2)'' are proved similarly.

We now establish $Z\in\mathcal{M}^2_T$ under (H1)'+(H2)'', the case (H1)+(H2)' being a specific instance. 
Applying It\^o's formula to $Y^{\eta}$ with $\eta>2\delta+1$, we obtain (after setting $t=0$ and rearranging)
\begin{align}
&\int_0^T\!\!\eta\left(\frac{\eta-1}{2}-\delta\right)Y_s^{\eta-2}|Z_s|^2ds \notag\\ &\leq \!-Y_0^{\eta}\!+\!X^{\eta}\!+\!\int_0^T\!\eta Y^{\eta-1}_s(\alpha_sY_s\!+\!\beta_sY_s|\ln(Y_s)|\!+\!\gamma_s|Z_s|)ds -\int_0^T\eta Y_s^{\eta-1}Z_sdW_s. 
\label{eq:Z2}
\end{align}
We first handle the regime where $Y$ is large. Taking $\eta\geq 2$ (so that $Y^{\eta-2}$ does not vanish), introducing a localizing sequence as in Proposition~\ref{prop:EUgen}, taking expectations, and applying Young's inequality with the splitting $Y^{\eta-1}|Z|=Y^{\eta/2}Y^{(\eta-2)/2}|Z|$, we obtain
\begin{align*}
&\mathbb{E}\Big[\int_0^{\tau_n}\eta\left(\frac{\eta-1}{2}-\delta\right)Y_s^{\eta-2}|Z_s|^2ds\Big]\\
&\leq\mathbb{E}\left[-Y_0^{\eta}+X^{\eta}\right] \\
&+\mathbb{E}\left[\int_0^T \left(\frac{(\eta\alpha_s)^q}{q}+\frac{Y_s^{p\eta}}{p}+\frac{(\eta\beta_s)^q}{q} +\frac{Y_s^{p\eta}|\ln^p(Y_s)|}{p}+\frac{(\eta\gamma_sY_s^{\eta/2})^2}{2\varepsilon}+\frac{\varepsilon}{2}Y_s^{\eta-2}|Z_s|^2 \right)ds\right].
\end{align*}
Choosing $\frac{\varepsilon}{2}<\eta(\frac{\eta-1}{2}-\delta)$ and setting $K:=\eta(\frac{\eta-1}{2}-\delta)-\frac{\varepsilon}{2}>0$, a further application of Young's inequality yields
\begin{align*}
K\mathbb{E}\left[\int_0^{\tau_n}Y_s^{\eta-2}|Z_s|^2ds\right]
&\leq 
K'\bigg(1+\mathbb{E}\bigg[\sup_{t\in[0,T]}Y_t^{\eta p+\varepsilon'}\bigg]\bigg),
\end{align*}
where we have used $y^m|\ln^m(y)|\leq K_{\varepsilon',m}+y^{m+\varepsilon'}$ for any $m>1$ and $\varepsilon'>0$.
Since we have $\mathbb{E}[\sup_{t\in[0,T]}Y_t^{p(2\delta+1)(e^B+1)}]<+\infty$ with $e^B+1\geq 2$, imposing $\eta p+\varepsilon'=2p(2\delta+1)$ requires $\eta<2(2\delta+1)$. Combined with $\eta\geq 2$ and $\eta>2\delta+1$, a valid $\eta$ exists for any $\delta\in\R_+$.
Thus, $\mathbb{E}[\int_0^TY_s^{\eta-2}|Z_s|^2ds]<+\infty$ by Fatou's lemma, and setting $A:=\{(\omega,t): Y_t(\omega)\geq C\}$ for $C>0$:
\begin{equation} 
C^{\eta-2}\mathbb{E}\left[\int_0^T|Z_s|^2\mathbb{I}_Ads\right]\leq\mathbb{E}\left[\int_0^TY_s^{\eta-2}|Z_s|^2ds\right]<+\infty.
\label{eq:yunbound}
\end{equation}

For the regime where $Y$ is small, we apply It\^o's formula to $\ln(1+Y_t)$. Since $g\geq 0$ and $Y>0$:
$$\ln(1+Y_t)\geq \ln(1+X)+\int_t^T\frac{1}{2(1+Y_s)^2}|Z_s|^2ds-\int_t^T\frac{1}{1+Y_s}Z_sdW_s.$$
Introducing the stopping time $\tau_n:=\inf\{t>0:\int_0^t\frac{|Z_s|^2}{(1+Y_s)^2}ds\geq n\}\wedge T$, taking expectations, using $\ln(1+x)\leq x$, and letting $n\to\infty$ via Fatou's lemma, we obtain 
$$\mathbb{E}\left[\int_0^T\frac{1}{2(1+Y_s)^2}|Z_s|^2ds\right]\leq Y_0+\mathbb{E}\left[X\right]<+\infty,$$
whence
\begin{equation} 
\frac{1}{2(1+C)^2}\mathbb{E}\left[\int_0^T|Z_s|^2\mathbb{I}_{A^c}ds\right]\leq Y_0+\mathbb{E}\left[X\right]<+\infty.
\label{eq:ybound}
\end{equation}
Combining~\eqref{eq:yunbound} and~\eqref{eq:ybound} gives $Z\in\mathcal{M}^2_T$.
\end{proof}
\begin{corollary}
Assume (H1)' with $h(t,y,z)=\alpha_t y+\beta_t y|\ln(y)|+\gamma_t|z|+\delta |z|^2/y$. 
If $\alpha,\beta,\gamma\in\mathcal{H}^{\infty}_T$ and there exists $\lambda\geq 2\|\beta\|^T_{\infty}$ such that $\mathbb{E}[X^{(2\delta+1)(e^{\lambda T}+1)}]<+\infty,$ then Equation~\eqref{eq:LNQ} admits a unique solution with regularity $(Y,Z)\in\mathcal{H}^{(2\delta+1)(e^{\lambda T}+1)}_T\times\mathcal{M}^2_T.$
\label{cor:EU1}
\end{corollary}
\begin{remark}
Proposition~\ref{prop:EUgen} and Theorem~\ref{th:Egen} also establish the existence of maximal and minimal solutions to Equation~\eqref{eq:LNQ} driven by $h$ and $g$, respectively, as follows from Lemma~\ref{lem:bahlali}.
\end{remark}

\subsection{Further regularities in the bounded case}
\begin{proposition}
Assume $X\in L^{\infty}_{+}(\mathcal{F}_T)$, $\alpha,\beta,\gamma\in \mathcal{H}^{\infty}_T$ and $\delta>0$. 
Let $h$ be as defined in Corollary~\ref{cor:EU1}.
Then the BSDE
\begin{equation}
Y^h_t=X+\int_t^Th(s,Y_s^h,Z_s^h)ds-\int_t^TZ^h_sdW_s,
\label{eq:exdr}
\end{equation}
admits a unique solution $(Y^h,Z^h)\in\mathcal{H}^{\infty}_T\times\text{BMO}(\mathbb{P})$. 
Furthermore, under the same hypotheses on the coefficients, if $g:[0,T]\times\Omega\times\R_+\times\R^n\to\R_+$ verifies (H1)', then the corresponding BSDE admits at least one solution $(Y,Z)\in \mathcal{H}^{\infty}_T\times\text{BMO}(\mathbb{P})$. 
Specifically, it admits maximal and minimal solutions with such regularity, among all possible solutions verifying $0<Y\leq Y^h$.
\label{prop:Linf}
\end{proposition}
\begin{corollary}
Consider $X\in L^{\infty}_{+}(\mathcal{F}_T)$, $\alpha,\beta\in\mathcal{H}^{\infty}_T$ and $\gamma\in\text{BMO}(\mathbb{P})$, and let $g:[0,T]\times\Omega\times\R_+\times\R^n\to\R_+$ verify (H1)''. 
Then, the BSDE with parameters $(X,g)$ admits at least one solution. 
Furthermore, any solution is such that $(Y,Z)\in\mathcal{H}^{\infty}_T\times\text{BMO}(\mathbb{P})$.
\label{cor:Linff}
\end{corollary}
\subsection{Uniqueness}
We establish uniqueness for Equation~\eqref{eq:LNQ} under the combined LN-Q growth rate, allowing for unbounded terminal conditions that are not necessarily bounded away from zero. 
Our approach, based on a comparison principle, extends \cite{BT19} (where only $|z|^2/y$-growth is treated under more restrictive assumptions, including bounded terminal conditions bounded away from zero) and complements \cite{BKK17} (where $y|\ln(y)|$-growth is handled via monotonicity conditions). 
We require the following assumptions:
\begin{itemize}
\item[A)]\label{item:A} $0\leq g(t,y,z) \leq \alpha_t y+\beta_t y|\ln(y)| +\gamma_t\cdot z+\delta |z|^2/y:=h(t,y,z), \ \ \forall(y,z)\in\R_+\times\R^n,$
    where $\alpha,\beta,\gamma\in\mathcal{H}^{\infty}_T$ and $\delta>0$.
\item[C)]\label{item:C} The driver $g$ is jointly convex in $(y,z)\in\R_+\times\R^n$.
\end{itemize}
One can equivalently write the bound in \hyperref[item:A]{A)} as
$0\leq \mathbb{I}_{\{y>0\}}\cdot g(t,|y|,z) \leq \mathbb{I}_{\{y>0\}}\cdot (\alpha_t y^++\beta_t y^+|\ln(y^+)| +\gamma_t\cdot z+\delta |z|^2/|y|)$ for all $(y,z)\in\R\times\R^n$,
with the convention that all members are $0$ when $y=0$.
The term ``$\gamma_t\cdot z$'' can also be replaced by ``$\gamma_t|z|$'' with $\gamma\in\mathcal{H}^{\infty}_T$ without altering any of the following results.

We start with a generalization of the stochastic Gr\"onwall inequality (see \cite{HKK22,WF18}) that accounts for the nonlinearity $y|\ln(y)|$.
\begin{lemma}
\label{lem:biharine}
Consider a positive process $\beta\in \mathcal{H}^{\infty}_T$ and $X\in L_+^{p(e^B+1)}(\mathcal{F}_T)$ for some $p>1$, where $B:=\|\beta\|^T_{\infty}\cdot T$. 
Let $u\in\mathcal{H}^{p(e^{B}+1)}_T$ such that $u_t\geq0 \ d\mathbb{P}\times dt$-a.s. 
If $u$ verifies $d\mathbb{P}$-a.s.\ for any $t\in[0,T]$:
$$u_t\leq\mathbb{E}\left[\left.X+\int_t^T\beta_su_s|\ln(u_s)|ds\right|\mathcal{F}_t\right],$$
then there exists an increasing function $\psi:\R_+\to\R$ such that $d\mathbb{P}$-a.s.\ for any $t\in[0,T]$:
$$u_t\leq\mathbb{E}\bigg[\psi^{-1}\left(\psi(X)+\int_t^T\beta_sds\right)\Bigg|\mathcal{F}_t\bigg],$$
where $\psi^{-1}:\text{Range}(\psi)\to\R_+$ is the inverse of $\psi$, given by
$$
\psi(x):=
\begin{cases}
\frac{x-2}{\ln(4)} \ &\text{ if } 0\leq x\leq2, \\
\ln(\ln(x))-\ln(\ln(2)) &\text{ if } x>2.
 \end{cases}
$$
\end{lemma}
\begin{theorem}
\label{th:EUg1}
Assume $g$ satisfies \hyperref[item:A]{A)}, \hyperref[item:C]{C)}, and $X\in L_+^{p(2\delta+1)(e^B+1)}$, where $p>\max\{\frac{2}{2\delta+1},1\}$ and $B=(\|\alpha\|^T_{\infty}+\|\beta\|^T_{\infty})T$. 
Then, the BSDE~\eqref{eq:LNQ} with parameters $(X,g)$ admits a unique solution $(Y,Z)\in\mathcal{H}_T^{p(2\delta+1)(e^B+1)}\times\mathcal{M}^2_T$ such that $0<Y\leq Y^h$, where $Y^h$ is the (maximal) solution corresponding to the driver $h$ as in \hyperref[item:A]{A)}.\footnote{Here, ``unique solution'' means that if $(Y,Z)$ and $(Y',Z')$ are two solutions to Equation~\eqref{eq:LNQ} with $0<Y,Y'\leq Y^h$, then $Y$ and $Y'$ are indistinguishable processes, and $Z_t=Z'_t$ $d\mathbb{P}\times dt$-a.s.}
\end{theorem}
The proof of Theorem~\ref{th:EUg1} relies on the following comparison principle.
\begin{proposition}
\label{prop:comparisongen}
With the same notation as in Theorem~\ref{th:EUg1}, consider a driver $g$ (resp.\ $g'$) verifying \hyperref[item:A]{A)}, \hyperref[item:C]{C)} and $X,X'\in L_+^{p(2\delta+1)(e^B+1)}$ with $X\leq X'$. 
If $(Y,Z), (Y',Z') \in\mathcal{H}_T^{p(2\delta+1)(e^B+1)}\times\mathcal{M}^2_T$ are solutions with $0<Y,Y'\leq Y^{h}$, and
$$g(t,Y'_t,Z'_t)\leq g'(t,Y'_t,Z'_t) \ \ \ (\text{resp. } g(t,Y_t,Z_t)\leq g'(t,Y_t,Z_t) ) \  d\mathbb{P}\times dt \text{-a.s.},$$
then $Y_t\leq Y'_t$ for any $t\in[0,T]$, $d\mathbb{P}$-a.s.
\end{proposition}
\begin{proof}
We treat the case $g(t,Y'_t,Z'_t)\leq g'(t,Y'_t,Z'_t)$; the other case is similar.
Fix $\theta\in(0,1)$ and define $P:=\frac{Y-\theta Y'}{1-\theta}$, $V:=\frac{Z-\theta Z'}{1-\theta}$. 
Then,
$$P_t=P_T+\int_t^TG(s,P_s,V_s)ds-\int_t^TV_sdW_s,$$
where, for $((1-\theta)y+\theta Y'_t(\omega),(1-\theta)z+\theta Z'_t(\omega))\in\R_+\times\R^n$,
\begin{align*}
G(\omega,t,y,z)\!:\!&=\!\frac{1}{1\!-\!\theta}\!\left[g(\omega,t,(1-\theta)y\!+\!\theta Y'_t(\omega),(1\!-\!\theta)z+\theta Z'_t(\omega))\!-\!\theta g(\omega,t,Y'_t(\omega),Z'_t(\omega))\right] \\
&\!+\!\frac{\theta}{1-\theta}\left[g(\omega,t,Y'_t(\omega),Z'_t(\omega))-g'(\omega,t,Y'_t(\omega),Z'_t(\omega))\right],
\end{align*}
and $G(\omega,t,y,z):=0$ otherwise.
By convexity of $g$, one verifies $\mathbb{I}_{\{y>0\}}G(t,y,z)\leq \mathbb{I}_{\{y>0\}}g(t,|y|,z)$ (with the convention that the right-hand side is zero when $y=0$). 
Applying It\^o-Tanaka's formula to $(P^+)^{\eta}$ with $\eta=\max\{2\delta+1,2\}$, using $\hyperref[item:A]{A)}$ and discarding the non-positive local time term, we obtain
\begin{align}
    &(P_t^+)^{\eta}\leq (P_T^+)^{\eta}+\int_t^T\eta(P_s^+)^{\eta}(\alpha_s+\beta_s|\ln(P_s^+)|)ds \notag\\
    &-\int_t^T\eta\mathbb{I}_{\{P_s>0\}}(P_s)^{\eta-1}V_sdW^{\gamma}_s,
\label{eq:condexp}
\end{align} 
where $W^{\gamma}$ denotes the Brownian motion under the Girsanov change of measure induced by $\gamma\in\mathcal{H}^{\infty}_T$.
Introducing the localization $\tau_n:=\inf\{s\geq t:\int_t^s(\eta\mathbb{I}_{\{P_u>0\}}(P_u^+)^{\eta-1}V_u)^2du\geq n\}\wedge T$, taking the $\mathbb{Q}^{\gamma}$-conditional expectation, and using $x^{\eta}\leq e^{\eta}+x^{\eta}|\ln(x)|$, we find
\begin{align*} 
(P_t^+)^{\eta}&\leq \mathbb{E}_{\mathbb{Q}^{\gamma}}\left[\left(\sup_{t\in[0,T]}P_{t}^+\right)^{\eta}+\int_0^{T}e^{\eta}\alpha_sds+\int_t^{\tau_n}(\alpha_s+\beta_s)(P_s^+)^{\eta}|\ln(P_s^+)^{\eta}|ds\bigg|\mathcal{F}_t\right].
\end{align*}
The conditional expectation under $\mathbb{Q}^{\gamma}$ is well-defined since $\mathcal{E}^{\gamma}_T\in L^m(\mathcal{F}_T)$ for any $m\geq 1$ (as $\gamma\in\mathcal{H}^{\infty}_T$), and the regularity of $Y,Y'$ ensures $(P^+)^{\eta}\in \mathcal{H}^{p'(e^B+1)}_T$ for some $p'>1$.
Applying Lemma~\ref{lem:biharine} (whose hypotheses can be verified via Young's inequality and the integrability of $P,\mathcal{E}^{\gamma}_T$), we get
\begin{align*}(P_t^+)^{\eta}\leq\mathbb{E}_{\mathbb{Q}^{\gamma}}\left[\psi^{-1}\left(\psi\left((P_{\tau_n}^+)^{\eta}+\int_0^{\tau_n}e^{\eta}\alpha_sds\right)+\int_t^{\tau_n}(\alpha_s+\beta_s)ds\right)\bigg|\mathcal{F}_t\right].
\end{align*}
Letting $n\to\infty$ via the dominated convergence theorem for conditional expectations, and using $(P_T^+)^{\eta}\leq|X|^{\eta}$ together with the monotonicity of $\psi,\psi^{-1}$ yields
$$\left(\left(Y_t-\theta Y_t'\right)^+\right)^{\eta}\leq(1-\theta)^{\eta}\mathbb{E}_{\mathbb{Q}^{\gamma}}\left[\psi^{-1}\left(\psi\left(|X|^{\eta}+\int_0^Te^{\eta}\alpha_sds\right)+\int_t^T(\alpha_s+\beta_s)ds\right)\bigg|\mathcal{F}_t\right].$$ 
Letting $\theta\to 1$ gives $Y_t\leq Y'_t$, $d\mathbb{P}$-a.s. 
\end{proof}
\begin{proof}[Proof of Theorem~\ref{th:EUg1}]
Existence follows from Theorem~\ref{th:Egen}. 
Uniqueness of $Y$ is immediate from Proposition~\ref{prop:comparisongen} and pathwise continuity. 
For uniqueness of $Z$, if $(Y,Z)$ and $(Y,Z')$ both solve~\eqref{eq:LNQ}, It\^o's formula applied to $(Y-Y)^2=0$ yields $\int_0^T|Z_s-Z'_s|^2ds=0$ $d\mathbb{P}$-a.s.
\end{proof}

\subsection{Stability}
The analysis of the stability of BSDE solutions has long been an important topic in the theory of BSDEs (e.g., \cite{EPQ97, K00, BH08}) and continues to attract attention in the recent literature, e.g., in \cite{PDS23}. 
We establish stability for the BSDE~\eqref{eq:LNQ} under assumptions~\hyperref[item:A]{A)} and~\hyperref[item:C]{C)}, building on the comparison principle of Proposition~\ref{prop:comparisongen}. 
This yields as a special case of interest the first stability result under $|z|^2/y$-growth.
We first derive refined regularity estimates for the $z$-component.
\begin{proposition}
With the same notation as in Theorem~\ref{th:Egen}, under assumption~\hyperref[item:A]{A)}, let $X\in L_+^{p(2\delta+1)(e^B+1)}$ with $p>1$ and $B=\|\beta\|^T_{\infty}\cdot T$. 
Any solution $(Y,Z)$ to~\eqref{eq:LNQ} with $0<Y\leq Y^h$ verifies, for some $C>0$ depending on $p,\|\alpha\|^T_{\infty},\|\beta\|^T_{\infty},\|\gamma\|^T_{\infty},\delta$ and $T$:
$$\mathbb{E}\left[\left(\int_0^T|Z_s|^2ds\right)^{p}\right]+\mathbb{E}\left[\left(\int_0^Tg(s,Y_s,Z_s)ds\right)^{2p}\right]\leq C\mathbb{E}\left[X^{p(2\delta+1)(e^{B}+1)}\right]<+\infty.$$ 
Furthermore, if $X\geq\varepsilon>0$ a.s., then
$$\mathbb{E}\left[\left(\int_0^T\frac{|Z_s|^2}{Y_s}ds\right)^{p}\right]\leq C'\mathbb{E}\left[X^{p(2\delta+1)(e^{B}+1)}\right]<+\infty,$$
for some $C'>0$ depending on $p,\|\alpha\|^T_{\infty},\|\beta\|^T_{\infty},\|\gamma\|^T_{\infty},\delta$ and $T$.
\label{prop:Zreg}
\end{proposition}
\begin{proof}
For brevity, we assume $\alpha,\beta,\gamma$ constant; the same results hold for bounded random coefficients as in~\hyperref[item:A]{A)}. 
Throughout, $C>0$ denotes a generic constant (depending on $T,\alpha,\beta,\gamma,\delta,p$) that may vary from line to line.

\emph{Step~1: $\mathcal{M}^{2p}_T$-regularity of $Z$.}
If $Y$ is bounded, then $Z\in\text{BMO}(\mathbb{P})$ by Proposition~\ref{prop:Linf}, hence $Z\in\mathcal{M}^p_T$ for any $p\geq 1$. 
For the general case, proceeding as in Theorem~\ref{th:Egen} with $\eta\geq 2$, $\eta>2\delta+1$, and applying Young's inequality, we obtain
\begin{align*}
K\int_0^TY_s^{\eta-2}|Z_s|^2ds\leq C\sup_{t\in[0,T]}Y_t^{\eta+\varepsilon'}-\int_0^T\eta Y_s^{\eta-1}Z_sdW_s,
\end{align*}
for suitable $K>0$ and $\varepsilon'>0$. 
Raising to the $p$-th power, taking expectations, and using BDG's inequality with the localization $\tau_n:=\inf\{t:\int_0^t(\eta Y_s^{\eta-1}Z_s)^2ds\geq n\}\wedge T$, we find
\begin{align}
K\mathbb{E}\left[\left(\int_0^{\tau_n}\eta Y_s^{\eta-2}|Z_s|^2ds\right)^p\right]&\leq C\mathbb{E}\left[\sup_{t\in[0,T]}Y_t^{p\eta+\varepsilon''}\right]+C\mathbb{E}\left[\left(\int_0^{\tau_n}\eta^2 Y_s^{2(\eta-1)}|Z_s|^2 ds\right)^{p/2}\right] \notag \\
&\leq C\mathbb{E}\left[\sup_{t\in[0,T]}Y_t^{p\eta+\varepsilon''}\right]+\frac{1}{2}\mathbb{E}\left[\left(\int_0^{\tau_n}Y_s^{\eta-2}|Z_s|^2ds\right)^{p}\right],\label{eq:Zpp}
\end{align}
where $\varepsilon''=p\varepsilon'$ and the last step uses Young's inequality with a suitable constant. 
Absorbing the last term into the left-hand side and letting $n\to\infty$ by Fatou's lemma:
$$\mathbb{E}\left[\left(\int_0^{T}Y_s^{\eta-2}|Z_s|^2ds\right)^p\right]\leq C\mathbb{E}\left[\sup_{t\in[0,T]}Y_t^{p\eta+\varepsilon''}\right].$$
Choosing $\eta,\varepsilon''$ as in Theorem~\ref{th:Egen} and using Proposition~3.2 of \cite{B20}:
$$\mathbb{E}\left[\left(\int_0^T|Z_s|^2ds\right)^p\right]\leq C\mathbb{E}\left[X^{p(2\delta+1)(e^{B}+1)}\right],$$
where the passage from $Y^{\eta-2}|Z|^2$ to $|Z|^2$ follows by the same splitting argument as in Theorem~\ref{th:Egen}.

\emph{Step~2: Regularity of the driver process.}
Rearranging~\eqref{eq:LNQ} and using $g\geq0$:
$$\left(\int_0^T|g(s,Y_s,Z_s)|ds\right)^{2p}\leq c_p\left(\sup_{t\in[0,T]}Y_t^{2p}+\Big|\int_0^TZ_sdW_s\Big|^{2p}\right).$$ 
Taking expectations and applying BDG's inequality:
$$\mathbb{E}\left[\left(\int_0^T|g(s,Y_s,Z_s)|ds\right)^{2p}\right]\leq c_p\mathbb{E}\left[\sup_{t\in[0,T]}Y_t^{2p}+\left(\int_0^T|Z_s|^2ds\right)^{p}\right]<+\infty.$$

\emph{Step~3.} If $X\geq\varepsilon>0$, then $Y_t\geq\varepsilon$ $d\mathbb{P}\times dt$-a.s., and the bound on $\int_0^T|Z_s|^2/Y_s\,ds$ follows immediately from Step~1.
\end{proof}
\begin{theorem}
\label{th:stability}
With the same notation as in Theorem~\ref{th:EUg1}, 
consider $(X^n,g^n)_{n\in\mathbb{N}}$ and $(X,g)$ with $g$ satisfying~\hyperref[item:A]{A)} and~\hyperref[item:C]{C)}. 
Let $\sup_ng^n$ verify~\hyperref[item:A]{A)}, $g^n$ verify~\hyperref[item:C]{C)} for any $n\in\mathbb{N}$, and $(X^n)_{n\in\mathbb{N}},X$ be strictly positive with $\sup_{n\in\mathbb{N}}X^n,X\in L^{p(2\delta+1)(e^{B}+1)}$, where $p>\max\{1,\frac{2}{2\delta+1}\}$. 
Denote by $(Y,Z)$ (resp.\ $(Y^n,Z^n)$) the unique solution to the BSDE with parameters $(X,g)$ (resp.\ $(X^n,g^n)$) such that $0< Y \leq Y^{h}$ (resp.\ $0< Y^n\leq Y^h$). 
Assume $X^n\to X$ $d\mathbb{P}$-a.s.\ and 
\begin{equation}\label{eq:driverconv}
\int_0^T|g^n(t,Y_t,Z_t)-g(t,Y_t,Z_t)|dt\to 0 \quad \text{in } L^{p(2\delta+1)(e^B+1)} \text{ as } n\to\infty.
\end{equation}
Then, $(Y^n,Z^n)\to (Y,Z)$ in $L^{q}\times \mathcal{M}^{2p}_T$ for any $q\in[1,{p(2\delta+1)(e^B+1)})$. 
\end{theorem}
\begin{proof}
By Theorem~\ref{th:Egen}, Proposition~\ref{prop:Zreg} and the assumptions on $(X^n,g^n)_{n\in\mathbb{N}}$:
\begin{equation*}
\sup_{n\in\mathbb{N}}\mathbb{E}\left[\sup_{t\in[0,T]}(Y_t^n)^{p(2\delta+1)(e^B+1)}+\left(\int_0^T|Z^n_s|^2 ds\right)^p\right]\leq\!K\sup_{n\in\mathbb{N}}\mathbb{E}\left[(X^{n})^{p(2\delta+1)(e^B+1)}\right]\!\!<\!+\infty,
\end{equation*}
for some $K>0$. 
By Vitali's convergence theorem, it suffices to show $\sup_{t\in[0,T]}|Y^n_t-Y_t|\to 0$ in probability.

We first show that convergence of $Y^n\to Y$ in $\mathcal{H}_T^{q}$ for any $q\in[1,p(2\delta+1)(e^B+1))$ implies $Z^n\to Z$ in $\mathcal{M}^{2p}_T$. 
By It\^o's formula and BDG's inequality:
    \begin{align*}
    &\mathbb{E}\left[\left(\int_0^T|Z_s^n-Z_s|^2ds\right)^p\right]\leq C\mathbb{E}\Bigg[|X^n-X|^{2p}+\sup_{t\in[0,T]}|Y^n_t-Y_t|^{2p}\Bigg] \\
    &+C\mathbb{E}\Bigg[\sup_{t\in[0,T]}|Y^n_t-Y_t|^{p}\left(\int_0^T|g^n(s,Y^n_s,Z^n_s)-g(s,Y_s,Z_s)|ds\right)^p\Bigg].
    \end{align*}
The first two terms converge to $0$ since $2p<p(2\delta+1)(e^B+1)$. 
For the third term, Proposition~\ref{prop:Zreg} and the integrability of $(X^n)_{n}$ give
$\sup_{n}\mathbb{E}[(\int_0^T|g^n(s,Y^n_s,Z^n_s)|ds)^{2p}]<+\infty,$
so H\"older's inequality yields
\begin{align*}
&\mathbb{E}\left[\sup_{t\in[0,T]}|Y^n_t-Y_t|^{p}\left(\int_0^T|g^n(s,Y^n_s,Z^n_s)-g(s,Y_s,Z_s)|ds\right)^p\right] \\
&\leq K\mathbb{E}\left[\sup_{t\in[0,T]}|Y^n_t-Y_t|^{2p}\right]^{\frac{1}{2}}\xrightarrow{n\to\infty} 0.
\end{align*}

It remains to prove $\sup_{t\in[0,T]}|Y^n_t-Y_t|\to 0$ in $L^{\eta}$ for $\eta=\max\{2,2\delta+1\}$.
We proceed as in Proposition~\ref{prop:comparisongen}, estimating $Y^n-\theta Y$. 
Define $P:=\frac{Y^n-\theta Y}{1-\theta}$, $V:=\frac{Z^n-\theta Z}{1-\theta}$, satisfying
\begin{equation}
P_t=P_T+\int_t^TG(s,P_s,V_s)ds-\int_t^TV_sdW_s,
\label{eq:stab1}
\end{equation}
where
\begin{align*}
G(\omega,t,y,z)\!:\!&=\!\frac{1}{1\!-\!\theta}\left[g^n(\omega,t,\!(1-\theta)y\!+\!\theta Y_t(\omega),\!(1-\theta)z\!+\!\theta Z_t(\omega))\!-\!\theta g^n(\omega,t,\!Y_t(\omega),Z_t(\omega))\right] \\
&+\frac{\theta}{1-\theta}\left[g^n(\omega,t,Y_t(\omega),Z_t(\omega))\!-\!g(\omega,t,Y_t(\omega),Z_t(\omega))\right],
\end{align*}
for $((1-\theta)y+\theta Y_t(\omega),(1-\theta)z+\theta Z_t(\omega))\in\R_+\times\R^n$, and $G:=0$ otherwise. 
By convexity of $g^n$, setting $\delta^n_{\theta} g(\fatdot{}):=\frac{\theta}{1-\theta}[g^n(\fatdot{},Y_{\fatdot{}},Z_{\fatdot{}})-g(\fatdot{},Y_{\fatdot{}},Z_{\fatdot{}})]$:
$$\mathbb{I}_{\{y>0\}}G(t,y,z)\leq \mathbb{I}_{\{y>0\}}\left(g^n(t,|y|,z)+|\delta^n_{\theta} g(t)|\right).$$
Following the same steps as in Proposition~\ref{prop:comparisongen}, we obtain the analog of~\eqref{eq:condexp}:
\begin{align*}
   (P_t^+)^{\eta} &\leq (P_T^+)^{\eta}+\sup_{t\in[0,T]}(P_t^+)^{\eta-1}\int_0^T \eta|\delta^n_{\theta}g(s)|ds \\
    &+\int_t^T \!\!\eta(P_s^+)^{\eta}(\alpha_s+\beta_s|\ln(P_s^+)|)ds -\int_t^T\eta\mathbb{I}_{\{P_s>0\}}(P_s)^{\eta-1}V_sdW^{\gamma}_s.
\end{align*}
Setting $\Gamma^n_T:=\sup_{t\in[0,T]}(P_t^+)^{\eta-1}\int_0^T|\delta^n_{\theta}g(s)|ds\leq C(\sup_{t}(P_t^+)^{\eta}+(\int_0^T|\delta^n_{\theta}g(s)|ds)^{\eta})$, which satisfies the hypotheses of Lemma~\ref{lem:biharine}, we apply this lemma and Doob's inequality (for $p'>1$ small enough) to obtain, assuming w.l.o.g.\ $\gamma\equiv0$:
\begin{align}
&\mathbb{E}\left[\left(\sup_{t\in[0,T]}(Y^n_t-Y_t)^+\right)^{\eta}\right] \notag \\ 
&\leq C_{p'}(1-\theta)^{\eta}\mathbb{E}\left[\left(P_T^+\right)^{p'\eta e^B}+\left(\Gamma_T^n\right)^{p'e^B}+\left(\int_0^Te^{\eta}\alpha_sds\right)^{p'e^B}\right].  \label{eq:stability}
\end{align}
For each fixed $\theta\in[0,1]$, the integrability of $X^n$ gives
$\mathbb{E}[(P_T^+)^{p'\eta e^B}]\to\mathbb{E}[X^{p'\eta e^B}]$ as $n\to\infty$.
Furthermore, by H\"older's inequality with conjugate exponents $l=\frac{\eta}{\eta-1}$ and $m=\eta$, using $L^{p'\eta e^B}$-boundedness of $\sup_t P_t^+$, we find
\begin{align*}
\mathbb{E}\left[(\Gamma^n_T)^{p'e^B}\right]
&\leq C\mathbb{E}\left[\left(\int_0^T|\delta^n_{\theta}g(s)|ds\right)^{p'\eta e^B}\right]^{\frac{1}{\eta}}\xrightarrow{n\to\infty} 0,
\end{align*}
by the convergence assumption~\eqref{eq:driverconv}. 
Letting $n\to\infty$ in~\eqref{eq:stability}:
\begin{equation*}
\limsup_{n\to\infty}\mathbb{E}\left[\sup_{t\in[0,T]}\left((Y^n_t-Y_t)^+\right)^{\eta}\right]\leq(1-\theta)^{\eta}\mathbb{E}\left[X^{p'\eta e^B}+\left(\int_0^Te^{\eta}\alpha_sds\right)^{p'e^B}\right].
\end{equation*}
Since the left-hand side does not depend on $\theta$, letting $\theta\to 1$ gives $\limsup_{n\to\infty}\mathbb{E}[\sup_{t\in[0,T]}((Y^n_t-Y_t)^+)^{\eta}]=0$. 

Similarly, estimating $Y-\theta Y^n$ with $P':=\frac{Y-\theta Y^n}{1-\theta}$, $V':=\frac{Z-\theta Z^n}{1-\theta}$ and $\delta^n_{\theta} g'(\fatdot{}):=\frac{1}{1-\theta}[g(\fatdot{},Y_{\fatdot{}},Z_{\fatdot{}})-g^n(\fatdot{},Y_{\fatdot{}},Z_{\fatdot{}})]$, the same argument yields
$\limsup_{n\to\infty}\mathbb{E}[\sup_{t\in[0,T]}((Y_t-Y^n_t)^+)^{\eta}]=0$. 
Combining both estimates and noting that the limit exists due to $\ds\liminf_{n\to\infty}\mathbb{E}[\sup_{t\in[0,T]}((Y_t-Y^n_t)^+)^{\eta}]\geq0$, we get $\lim_{n\to\infty}\mathbb{E}[\sup_{t\in[0,T]}|Y_t-Y^n_t|^{\eta}]=0$.
\end{proof}
The condition~\eqref{eq:driverconv} is automatically satisfied when $g^n\equiv g$ for all $n$. 
The following corollary provides further sufficient conditions.
\begin{corollary}
With the same notation and hypotheses as in Theorem~\ref{th:stability}, let $X^n\to X$ $d\mathbb{P}$-a.s.\ and $\sup_{n}X^n,X\in L^p$ for any $p\geq 1.$ 
Assume that, $d\mathbb{P}\times dt$-a.s.\ for any $(y,z)\in\R_+\times\R^n$, $g^n(t,y,z)\to g(t,y,z)$, and that there exists $\varepsilon>0$ such that $X\geq \varepsilon$ $d\mathbb{P}$-a.s. 
Then, $(Y^n,Z^n)\to (Y,Z)$ in $\mathcal{H}^p_T\times \mathcal{M}^p_T$ for any $p\geq 1.$
\label{cor:stability}
\end{corollary}
\begin{remark}
The assumption $\sup_nX_n,X\in L^p$ for any $p\geq1$ in Corollary~\ref{cor:stability} is natural, since \cite{BH08} assumed the finiteness of all exponential moments for quadratic BSDEs (see \cite[Proposition~7]{BH08}), and Proposition~\ref{prop:equivalence} suggests that exponential regularity for quadratic BSDEs leads to $L^p$ regularity under LN-Q growth.
The hypothesis $X\geq\varepsilon$ is also natural when seeking further regularities, as shown in \cite{BT19} and, within the framework of return risk measures, in \cite{BLR18,LR22}. 
\end{remark}
\setcounter{equation}{0}


\section{Two-Driver BSDEs}\label{sec:2DBSDEs}

In this section, we develop a comprehensive framework for two-driver BSDEs of the form
\begin{equation}
    Y_t=X+\int_t^Tg_1(s,Y_s,Z_s)ds-\int_t^Tg_2(s,Y_s,Z_s)dW_s,
    \label{eq:twodrivers}
\end{equation}
which embeds the GBSDE in Equation~\eqref{eq:GBSDE} as a special case. 
Here, $g_1$ satisfies the growth rate 
\begin{equation} 
  0\leq g_1(t,y,z)\leq y(\alpha_t+\beta_t|\ln(y)|+\gamma_t|z|+\delta|z|^2)=:h_1(t,y,z), \ \ \forall (y,z)\in\R_+\times\R^n,
  \label{eq:GGR}
\end{equation}
with $\alpha,\beta,\gamma$ non-negative predictable processes and $\delta>0$, and $g_2:[0,T]\times\Omega\times\R_+\times\R^n\to\R^n$ is $\mathcal{P}\times\mathcal{B}(\R_+)\times\mathcal{B}(\R^n)/\mathcal{B}(\R^n)$-measurable.
A couple $(Y,Z)$ is a \emph{solution} to Equation~\eqref{eq:twodrivers} if it satisfies the equation in the It\^o sense, $Y$ is continuous and predictable, $Z\in\mathcal{L}^2_T$ is predictable, and $\int_0^Tg_1(s,Y_s,Z_s)ds<+\infty$, $\int_0^T|g_2(s,Y_s,Z_s)|^2ds<+\infty$ hold $d\mathbb{P}$-a.s.

The following proposition shows that, under a suitable growth condition on $g_2$, Equation~\eqref{eq:twodrivers} can be reduced to an ordinary BSDE with LN-Q growth rate.

\begin{proposition}\label{prop:transBSDE}
Let $g_1$ be $\mathcal{P}\times\mathcal{B}(\R_+)\times\mathcal{B}(\R^n)/\mathcal{B}(\R_+)$-measurable and satisfy~\eqref{eq:GGR}, and let $g_2$ be as above. 
If there exists $K>0$ such that $|g_2(t,y,z)| \geq Ky|z|$ $d\mathbb{P}\times dt$-a.s.\ for all $(y,z)\in\R_+\times\R^n$, with $g_2$ injective and continuous in $z$, then Equation~\eqref{eq:twodrivers} admits a positive solution whenever the ordinary BSDE~\eqref{eq:LNQ} with LN-Q growth rate does.
\end{proposition}
\begin{proof}
The bound $|g_2(t,y,z)|\geq Ky|z|$ together with the injective continuity of $g_2$ in $z$ ensures that $z\mapsto g_2(\omega,t,y,z)$ is a bijection of $\R^n$ onto itself (cf.\ \cite[Theorem~59(V)]{Pr05}). 
We denote its $z$-inverse by $g_2^{-1}(\omega,t,y,\cdot)$, which satisfies
\begin{equation}\label{eq:invbound}
|g_2^{-1}(t,y,v)|\leq \frac{|v|}{Ky}, \quad \forall (y,v)\in\R_+\times\R^n, \ d\mathbb{P}\times dt\text{-a.s.}
\end{equation}
Substituting $V_{\fatdot{}}:=g_2(\fatdot{},Y_{\fatdot{}},Z_{\fatdot{}})$, Equation~\eqref{eq:twodrivers} becomes
\begin{equation}
Y_t=X+\int_t^Tg(s,Y_s,V_s)ds-\int_t^TV_sdW_s,
\label{eq:LNQin}
\end{equation}
where $g(t,y,v):=g_1(t,y,g_2^{-1}(t,y,v))$. Using~\eqref{eq:GGR} and~\eqref{eq:invbound}, one verifies
$$0\leq g(t,y,v)\leq \alpha_t y+\beta_ty|\ln(y)|+\frac{\gamma_t|v|}{K}+\frac{\delta|v|^2}{K^2y},$$
so that $g$ has the LN-Q growth rate. 
The $\mathcal{P}\times\mathcal{B}(\R_+)\times\mathcal{B}(\R^n)/\mathcal{B}(\R_+)$-measurability of $g$ follows from the measurability of $(\omega,t,y,v)\mapsto(\omega,t,y,g_2^{-1}(\omega,t,y,v))$, which is established by adapting \cite[Theorem~4.1]{PP90}: taking w.l.o.g.\ $\Omega=C^0([0,T],\R^n)$ with canonical Brownian motion, the map $G:(\omega,t,y,z)\mapsto(\omega,t,y,g_2(\omega,t,y,z))$ is a Borel bijection on a complete separable metric space, whence its inverse is measurable by \cite[Theorem~10.5]{EK86}.
Thus, if $(Y,V)$ solves~\eqref{eq:LNQin}, then $(Y,g_2^{-1}(\fatdot{},Y,V))$ solves~\eqref{eq:twodrivers}. 
We note, however, that the regularity of $Z_{\fatdot{}}=g_2^{-1}(\fatdot{},Y_{\fatdot{}},V_{\fatdot{}})$ depends on the properties of $g_2$ and must be analyzed separately.
\end{proof}

Informed by the reduction above, we impose the following standing assumptions.
\begin{itemize}
    \item[(G1)] $g_1:[0,T]\times\Omega\times\R_+\times\R^n\to\R_+$ is $\mathcal{P}\times\mathcal{B}(\R_+)\times\mathcal{B}(\R^n)/\mathcal{B}(\R_+)$-measurable, satisfies~\eqref{eq:GGR}, and is continuous in $(y,z)$.
    \item[(G2)] $g_2:[0,T]\times\Omega\times\R_+\times\R^n\to\R^n$ is $\mathcal{P}\times\mathcal{B}(\R_+)\times\mathcal{B}(\R^n)/\mathcal{B}(\R^n)$-measurable with $|g_2(t,y,z)| \geq Ky|z|=:h_2(t,y,z)$ for some $K>0$. 
    Moreover, $g_2$ is injective and continuous in $z$, and its $z$-inverse $g^{-1}_2$ (cf.\ Proposition~\ref{prop:transBSDE}) is continuous in $(y,z)$, $d\mathbb{P}\times dt$-a.s.
    \item[(G3)] Let $p>1$, let $\alpha,\beta,\gamma,\delta,K$ be as in~\eqref{eq:GGR} and (G2), and set $k\in(0,\frac{1\wedge(p-1)}{2})$, $A:=\int_0^T\alpha_tdt$, $B:=\int_0^T\beta_tdt$. Then
    \begin{equation*}
       \mathbb{E}\left[(1+X^{\frac{2\delta}{K^2}+1})^{p(e^{B}+1)}\exp\left(p(e^{B}+1)\left((A+B)+\frac{1}{4k}\int_0^T(\gamma_t/K)^2dt\right)\right)\right]<+\infty.
    \end{equation*}
 Furthermore, there exists $q'>0$ such that $\mathbb{E}[\int_0^Te^{q'\gamma_t}dt]<+\infty$ and $\alpha,\beta\in \mathcal{H}^{q}_T$ with $\frac{1}{q}+\frac{1}{p}=1$.
\end{itemize}

\begin{proposition}\label{prop:twodriverE}
Under (G1)--(G3), Equation~\eqref{eq:twodrivers} admits at least one solution $(Y,Z)$ with $0<Y\leq Y^{h_1,h_2}$, where $Y^{h_1,h_2}$ solves~\eqref{eq:twodrivers} with parameters $(X,h_1,h_2)$. Moreover,
$$\mathbb{E}\left[\ds\sup_{t\in[0,T]}Y_t^{p\left(\frac{2\delta}{K^2}+1\right)(e^{B}+1)}\right]<+\infty \quad \text{and} \quad Z\in\mathcal{L}^2_T.$$
If additionally $X\geq\varepsilon$ for some $\varepsilon>0$, then $Z\in\mathcal{M}^2_T$. 
\end{proposition}
\begin{proof}
By Proposition~\ref{prop:transBSDE}, it is sufficient to solve the ordinary BSDE~\eqref{eq:LNQin} with $g(t,y,v):=g_1(t,y,g_2^{-1}(t,y,v))$, which is continuous in $(y,v)$ and verifies $0\leq g(t,y,v)\leq \alpha_ty+\beta_ty|\ln(y)|+\frac{\gamma_t}{K}|v|+\frac{\delta}{K^2}\frac{|v|^2}{y}$. 
Theorem~\ref{th:Egen} yields a solution $(Y,V)$ with the stated regularity for $Y$, $V\in\mathcal{M}^2_T$, and $0<Y\leq Y^h$ where $h(t,y,v):=h_1(t,y,h_2^{-1}(t,y,v))$.
Setting $Z_{\fatdot{}}:=g_2^{-1}(\fatdot{},Y_{\fatdot{}},V_{\fatdot{}})$ gives a solution to~\eqref{eq:twodrivers} with $Z\in\mathcal{L}^2_T$.

It remains to show $Z\in\mathcal{M}^2_T$ when $X\geq\varepsilon>0$. 
First, $Y_t\geq\varepsilon$ for all $t\in[0,T]$, $d\mathbb{P}$-a.s., since after a suitable localization $Y_t=\mathbb{E}[X+\int_t^{\tau_n}g_1(s,Y_s,Z_s)ds|\mathcal{F}_t]\geq\mathbb{E}[X|\mathcal{F}_t]\geq\varepsilon$ by the positivity of $g_1$.
Next, applying It\^o's formula to $Y^{\eta}$ (valid since $Y>0$) with $\eta>\frac{2\delta}{K^2}+1$, and using $|g_2(t,y,z)|\geq Ky|z|$, we obtain
\begin{align}
&\int_0^T\eta\left(\frac{K^2(\eta-1)}{2}-\delta\right)Y_s^{\eta}|Z_s|^2ds \notag\\
&\leq -Y_0^{\eta}+X^{\eta}+\int_0^T\eta Y^{\eta}_s(\alpha_s+\beta_s|\ln(Y_s)|+\gamma_s|Z_s|)ds -\int_0^T\eta Y_s^{\eta-1}g_2(s,Y_s,Z_s)dW_s. \label{eq:Z22}
\end{align}
Since $Y\geq\varepsilon$, setting $K':=\varepsilon^{\eta}>0$, localizing, taking expectations, and applying Young's inequality yields
\begin{align*}
K''\,\mathbb{E}\left[\int_0^{\tau_n}|Z_s|^2ds\right]
&\leq C\left(1+\mathbb{E}\left[\sup_{t\in[0,T]}Y_t^{\eta p+\varepsilon'}+\sup_{t\in[0,T]}Y_t^{2\eta l}\right]\right),
\end{align*}
for $K''>0$, $1<l<p$ with $\frac{1}{m}+\frac{1}{l}=1$, and a suitable $\varepsilon'>0$.
Since $\mathbb{E}[\sup_{t}Y_t^{p(\frac{2\delta}{K^2}+1)(e^B+1)}]<+\infty$ and $e^B+1\geq 2$, the right-hand side is finite for $\eta$ satisfying: $\frac{2\delta}{K^2}+1<\eta<2(\frac{2\delta}{K^2}+1)$ and $\eta\leq\frac{p}{l}(\frac{2\delta}{K^2}+1)$. 
Such $\eta$ exists since $\frac{p}{l}>1$. Fatou's lemma gives $\mathbb{E}[\int_0^T|Z_s|^2ds]<+\infty$.
\end{proof}

\begin{remark}
If $\gamma\equiv0$ and (G3) is replaced by the weaker condition 
$\mathbb{E}[(1+X^{\frac{2\delta}{K^2}+1})^{e^B}\exp(e^{B}(A+B))]<+\infty,$
a positive solution still exists with $0<\sup_{t}\mathbb{E}[Y_t^{(2\delta+1)e^B/K^2}]<+\infty$ and $Z\in\mathcal{L}^2_{T}$, by Proposition~\ref{prop:EUgen}~(i).
\end{remark}

We now turn to uniqueness. By Proposition~\ref{prop:transBSDE} and Theorem~\ref{th:EUg1}, suitable hypotheses on the composition $g_1\circ g_2^{-1}$ ensure uniqueness for~\eqref{eq:twodrivers}, covering in particular Equation~\eqref{eq:GBSDE} with unbounded terminal conditions and without a monotonicity condition on $\tilde f$ w.r.t.\ $y$. 
We require:
\begin{itemize}
    \item[A')]\label{item:A'} $g_1$ is $\mathcal{P}\times\mathcal{B}(\R_+)\times\mathcal{B}(\R^n)/\mathcal{B}(\R_+)$-measurable with $$0\leq g_1(t,y,z) \leq \alpha_t y+\beta_t y|\ln(y)|+\delta y|z|^2, \ \ \forall(y,z)\in\R_+\times\R^n,$$
    where $\alpha,\beta$ are bounded non-negative stochastic processes and $\delta>0$.
    \item[C')]\label{item:C'} Given $g_1$ as in \hyperref[item:A']{A')} and $g_2$ as in (G2), the function $(\omega,t,y,v)\mapsto g_1(\omega,t,y,g_2^{-1}(t,\omega,y,v))$ is convex in $(y,v)$, $d\mathbb{P}\times dt$-a.s.
\end{itemize}
One can also include the term ``$y\gamma_t\cdot z$'' in \hyperref[item:A']{A')} with $\gamma\in\mathcal{H}^{\infty}_T$ without altering the results below. 

\begin{proposition}\label{prop:2DU}
Under \hyperref[item:A']{A')} and \hyperref[item:C']{C')}, let $X\in L_+^{p(\frac{2\delta}{K^2}+1)(e^B+1)}$ with $p>\max\{\frac{2K^2}{2\delta+K^2},1\}$ and $B=(\|\alpha\|^T_{\infty}+\|\beta\|^T_{\infty})T$. 
Then~\eqref{eq:twodrivers} admits a unique solution $(Y,Z)\in\mathcal{H}_T^{p(\frac{2\delta}{K^2}+1)(e^B+1)}\times\mathcal{L}^2_{T}$ among solutions with $0<Y\leq Y^{h_1,h_2}$.\footnote{Uniqueness means: if $(Y,Z)$ and $(Y',Z')$ both solve~\eqref{eq:twodrivers} with $0<Y,Y'\leq Y^{h_1,h_2}$, then $Y$ and $Y'$ are indistinguishable and $Z=Z'$ $d\mathbb{P}\times dt$-a.s.}
\end{proposition}
The following comparison principle is proved by analogous arguments and its proof is omitted for brevity.
\begin{proposition}\label{prop:2Dcomparisongen}
With the notation of Proposition~\ref{prop:2DU}, consider drivers $g_1,g_2$ (resp.\ $g'_1,g'_2$) satisfying \hyperref[item:A']{A')}, \hyperref[item:C']{C')} and $X,X'\in L_+^{p(2\delta+1)(e^B+1)}$ with $X\leq X'$. 
If $(Y,Z)\in\mathcal{H}_T^{p(2\delta+1)(e^B+1)}\times\mathcal{L}^2_T$ solves the BSDE with parameters $(X,g_1,g_2)$, $(Y',Z')\in\mathcal{H}_T^{p(2\delta+1)(e^B+1)}\times\mathcal{L}^2_T$ solves the BSDE with parameters $(X',g'_1,g'_2)$, $0<Y,Y'\leq Y^{h_1,h_2}$, and
\begin{align*}
g_1(t,Y'_t,g_2^{-1}(t,Y'_t,Z'_t))\leq g'_1(t,Y'_t,(g'_2)^{-1}(t,Y'_t,Z'_t)) \  d\mathbb{P}\times dt \text{-a.s.}
\end{align*}
(resp.\ with $(Y,Z)$ in place of $(Y',Z')$), then $Y_t\leq Y'_t$ for all $t\in[0,T]$, $d\mathbb{P}$-a.s.
\end{proposition}

\begin{corollary}\label{cor:euGBSDE}
(i) Consider Equation~\eqref{eq:GBSDE} with $0\leq\tilde f(t,y,z)\leq \alpha_t+\beta_t|\ln(y)|+\delta|z|^2=:\tilde h(t,y,z)$, where $\alpha,\beta$ and $X$ verify (G3) with $K=1$ and $\gamma\equiv0$. 
Then there exists at least one solution $(Y,Z)$ with $0<Y\leq Y^{\tilde h}$ and
$$\mathbb{E}\bigg[\ds\sup_{t\in[0,T]}Y_t^{p(2\delta+1)(e^{B}+1)}\bigg]<+\infty, \quad Z\in\mathcal{L}^2_{T}.$$
If $X\geq\varepsilon>0$, then $Z\in\mathcal{M}^2_T$. 

\noindent (ii) If additionally $\tilde f(t,y,z)=y(\tilde f_1(t,y)+\tilde f_2(t,z))$ with $\tilde f_1,\tilde f_2>0$, $\tilde f_1(\fatdot{},y)$ convex in $y$, $\tilde f_2(\fatdot{},z)$ convex in $z$, $\alpha,\beta\in\mathcal{H}^{\infty}_T$, and $p\geq\max\{\frac{2}{2\delta+1},1\}$, then the solution is unique among those with $0<Y\leq Y^{\tilde h}$. 
In particular, when $\tilde f$ is independent of $y$ and convex in $z$, uniqueness holds without further assumptions.
\end{corollary}

We conclude with a stability result for two-driver BSDEs.
\begin{proposition}\label{prop:2Dstability}
With the notation of Proposition~\ref{prop:2DU}, let $(X^n,g_1^n,g_2^n)_{n\in\mathbb{N}}$ and $(X,g_1,g_2)$ be such that $g_1$, $\sup_n g_1^n$ satisfy \hyperref[item:A']{A')}, $g_2$, $\sup_n g_2^n$ satisfy (G2), and $(g_1^n,g_2^n)_n$, $(g_1,g_2)$ satisfy \hyperref[item:C']{C')}. 
Assume $X^n,X\geq\varepsilon>0$, $X,\sup_n X^n\in L^{p(\frac{2\delta}{K^2}+1)(e^B+1)}$ with $p>\max\{\frac{2K^2}{2\delta+K^2},1\}$ and $B=(\|\alpha\|^T_{\infty}+\|\beta\|^T_{\infty})T$. 
Denote by $(Y,Z)$ (resp.\ $(Y^n,Z^n)$) the unique solution to~\eqref{eq:twodrivers} with parameters $(X,g_1,g_2)$ (resp.\ $(X^n,g_1^n,g_2^n)$) satisfying $0\leq Y^n,Y\leq Y^{h_1,h_2}$. 
Suppose that for some $C>0$, all $(g_2^n)^{-1}$ satisfy $|(g^n_2)^{-1}(t,y_1,z_1)-(g^n_2)^{-1}(t,y_2,z_2)|\leq C|z_1/y_1-z_2/y_2|$ $d\mathbb{P}\times dt$-a.s., and $(g_2^n)^{-1}(t,y,z)\to g_2^{-1}(t,y,z)$ pointwise $d\mathbb{P}\times dt$-a.s.

Setting $V_{\fatdot{}}:=g_2(\fatdot{},Y_{\fatdot{}},Z_{\fatdot{}})$, if $\int_0^Tg_1^n(t,Y_t,(g_2^n)^{-1}(t,Y_t,V_t))dt\to\int_0^Tg_1(t,Y_t,g_2^{-1}(t,Y_t,V_t))dt$ in $L^{p(\frac{2\delta}{K^2}+1)(e^B+1)}$ and $X_n\to X$ $d\mathbb{P}$-a.s., then $(Y^n,Z^n)\to(Y,Z)$ in $L^q\times \mathcal{M}^{r}_T$ for all $q\in[1,p(\frac{2\delta}{K^2}+1)(e^B+1))$ and $r\in[1,2p)$. 
\end{proposition}
\begin{proof}
Define $g^n(t,y,v):=g_1^n(t,y,(g_2^n)^{-1}(t,y,v))$ and $g(t,y,v):=g_1(t,y,g_2^{-1}(t,y,v))$. 
The ordinary BSDEs with parameters $(X^n,g^n)_n$ and $(X,g)$ satisfy all hypotheses of Theorem~\ref{th:stability}, which yields $(Y^n,V^n)\to(Y,V)$ in $L^q\times\mathcal{M}^{2p}_T$.
For the convergence of $Z^n=( g_2^n)^{-1}(\fatdot{},Y^n,V^n)$ to $Z=g_2^{-1}(\fatdot{},Y,V)$, the triangle inequality gives
\begin{equation}
|Z^n_t-Z_t|\leq|(g_2^n)^{-1}(t,Y^n_t,V^n_t)-(g_2^{n})^{-1}(t,Y_t,V_t)|+|(g_2^{n})^{-1}(t,Y_t,V_t)-g_2^{-1}(t,Y_t,V_t)|.
\label{eq:conveprob}
\end{equation}
The Lipschitz condition on $(g_2^n)^{-1}$, the bound $Y^n,Y\geq\varepsilon$, and the convergences $Y^n\to Y$, $V^n\to V$ in probability imply that the first term vanishes in probability. 
The second term converges to $0$ $d\mathbb{P}\times dt$-a.s.\ by the pointwise convergence of $(g_2^n)^{-1}$ to $g_2^{-1}$. 
Hence $Z^n\to Z$ in probability. 
Furthermore, since $Y^n\geq\varepsilon$ and $|(g_2^n)^{-1}(t,y,v)|\leq |v|/(Ky)$, we have
$$\sup_{n}\mathbb{E}\left[\left(\int_0^T|Z^n_t|^2dt\right)^{p}\right]\leq K_{\varepsilon}\sup_{n}\mathbb{E}\left[\left(\int_0^T|V^n_t|^2dt\right)^{p}\right]<+\infty,$$
so that $(Z^n)_n$ is uniformly integrable in $\mathcal{M}^r_T$ for $r\in[1,2p)$. Vitali's convergence theorem concludes.
\end{proof}

\begin{remark}\label{rem:stabless}
If $(X^n)_n,X$ are strictly positive but not necessarily bounded away from $0$, the convergence $Y^n\to Y$ in $\mathcal{H}^q_T$ for $q\in[1,p(2\delta+1)(e^B+1))$ is preserved, whereas the regularities for $(Z^n)_n$ and $Z$ do not necessarily hold.
\end{remark}
\setcounter{equation}{0}
\section{A Stochastic Control Application}\label{sec:app_lnq_control}

In this section, we show that the two-driver BSDE framework of Sections~\ref{sec:gbsdeLNQ}--\ref{sec:2DBSDEs} arises endogenously from a log-type stochastic differential utility (SDU) control problem \`a la Duffie--Epstein (\cite{DE92,LQ03}).
\subsection{Formulation}

Let \(W\) be an \(n\)-dimensional Brownian motion on the stochastic basis of Section~\ref{sec:prel}, let \(\theta\in\mathcal{H}^{\infty}_T\) be a bounded $\R^n$-valued predictable market price of risk, and set
\[
\A:=\Big\{p:\Omega\times[0,T]\to\R^n \text{ predictable}:\ 
f_t^p:=p_t\cdot\theta_t-\tfrac12|p_t|^2\ge0
\;\; d\mathbb{P}\times dt\text{-a.s.}\Big\}.
\]
For \(p\in\A\) and initial endowment \(x>0\), the wealth process \(X^p\) satisfies
\[
X_t^p=x + \int_0^t X_s^p\,p_s\cdot(\theta_s\,ds+dW_s),
\]
so that \(X^p>0\) and
\begin{equation}\label{eq:sec7_logwealth_new}
d\ln(X_t^p)=p_t\cdot dW_t+f_t^p\,dt.
\end{equation}
Since \(f_t^p=-\frac12|p_t-\theta_t|^2+\frac12|\theta_t|^2\), the feedback control \(p=\theta\) belongs to \(\A\), and every admissible control satisfies \(|p_t|\le 2|\theta_t|\) and \(0\le f_t^p\le \frac12|\theta_t|^2\), \(d\mathbb{P}\times dt\)-a.s.

Let \(\xi\in L_+^q(\mathcal{F}_T)\) for every \(q\ge1\). Furthermore, let \(\Lambda:[0,T]\times\Omega\to\R^{n\times n}\) be a predictable process such that, for some constants \(0<\underline\lambda\le\overline\lambda<+\infty\),
\begin{equation}
\label{eq:sec7_lambda_bounds}
\underline\lambda\,|z|
\le |\Lambda_t z|
\le \overline\lambda\,|z|,
\qquad
\forall z\in\R^n,
\quad d\mathbb{P}\times dt\text{-a.s.}
\end{equation}
As we will see, the matrix \(\Lambda\) describes how Brownian risk is transmitted to the SDU's opportunity process: the continuation-value sensitivity \(Z^p\) and the effective diffusion coefficient \(\Lambda Z^p\) of the utility process need not coincide.
The recursive utility/SDU performance criterion associated with \(p\in\A\) is defined by the following BSDE, consisting as usual of the terminal utility condition, aggregator term and diffusion term:
\begin{equation}\label{eq:sec7_primitive_bsde_new}
V_t^p =\ln(X_T^p\xi)+\int_t^T\Big(\ln\big(1+e^{V_s^p-\ln(X_s^p)}\big)+\big|\Lambda_s Z_s^p-p_s\big|^2\Big)ds-\int_t^T\Lambda_s Z_s^p\,dW_s.
\end{equation}
Here, \(\Lambda Z^p\) is the Brownian diffusion coefficient of the continuation utility \(V^p\), and \(|\Lambda Z^p-p|^2\) measures the residual continuation-utility variance after removing the direct wealth exposure \(p\) (see \cite{DE92}).
For each fixed \(p\in\A\), let \(V^p(X_T^p\xi)\) denote the first component, when it exists, of the solution to \eqref{eq:sec7_primitive_bsde_new} with logarithmic terminal condition \(\ln(X_T^p\xi)\). The control problem is then given by
\begin{equation}\label{eq:sec7_value_new}
\nu(x,\xi):=\sup_{p\in\A}V_0^p(X_T^p\xi).
\end{equation}
\subsection{Verification and optimal control}

For any solution $(V^p,\widetilde Z^p)$ of~\eqref{eq:sec7_primitive_bsde_new}, define the opportunity process
$Y_t^p:=e^{V_t^p}/X_t^p>0,$
which captures the continuation value per unit of current wealth.
The following theorem shows that $Y^p$ solves a two-driver BSDE, that this equation is well-posed by the theory of Section~\ref{sec:2DBSDEs}, and that $p^*=\theta$ is optimal.

\begin{theorem}[Verification and optimal control]\label{th:sec7_verification_new}
Under~\eqref{eq:sec7_lambda_bounds} and $\xi\in L_+^q(\mathcal{F}_T)$ for all $q\ge1$, the following statements hold.
\begin{enumerate}
    \item For every $p\in\A$, equation~\eqref{eq:sec7_primitive_bsde_new} admits a unique solution $(V^p,\widetilde Z^p)$ in the sense of Definition~\ref{def:solBSDE}, with the factorization
    \begin{equation}\label{eq:sec7_value_factorization_new}
    V_t^p=\ln(X_t^pY_t^p),\qquad \widetilde Z_t^p=p_t+\Lambda_t Z_t^p,\qquad t\in[0,T],
    \end{equation}
    where $(Y^p,Z^p)$ is the unique solution to the two-driver BSDE
    \begin{equation}\label{eq:sec7_twodriver_new}
    Y_t^p=\xi+\int_t^T g_1^p(s,Y_s^p,Z_s^p)\,ds-\int_t^T g_2(s,Y_s^p,Z_s^p)\,dW_s,
    \end{equation}
    with
    \begin{equation}\label{eq:sec7_g1g2_new}
    g_1^p(t,y,z):=y\ln(1+y)+f_t^p\,y+\tfrac12\,y\,|\Lambda_t z|^2,\qquad g_2(t,y,z):=y\Lambda_t z.
    \end{equation}

    \item The control $p^\star=\theta$ is optimal for~\eqref{eq:sec7_value_new}:
    \begin{equation}\label{eq:sec7_value_formula_new}
    \nu(x,\xi)=V_0^\theta(X_T^\theta\xi)=\ln(xY_0^\theta).
    \end{equation}
\end{enumerate}
\end{theorem}

\begin{proof}
\textit{Step~1: Reduction to a two-driver BSDE.}
Fix $p\in\A$ and write $X:=X^p$, $f_t:=f_t^p$.
Setting $R_t:=V_t-\ln(X_t)$, $Y_t:=e^{R_t}$, $\Gamma_t:=\widetilde Z_t-p_t$, and $Z_t:=\Lambda_t^{-1}\Gamma_t$, the dynamics of~\eqref{eq:sec7_primitive_bsde_new} combined with~\eqref{eq:sec7_logwealth_new} give
$dR_t=-(\ln(1+e^{R_t})+f_t+|\Gamma_t|^2)\,dt+\Gamma_t\,dW_t.$
Applying It\^o's formula to $Y_t=e^{R_t}$ yields
\[
dY_t=-\Big(Y_t\ln(1+Y_t)+f_tY_t+\tfrac12 Y_t|\Gamma_t|^2\Big)dt+Y_t\Gamma_t\,dW_t.
\]
Since $\Gamma_t=\Lambda_t Z_t$, this is exactly~\eqref{eq:sec7_twodriver_new}--\eqref{eq:sec7_g1g2_new}. 
Thus every solution of the primitive equation induces a solution of the two-driver BSDE, and conversely.

\smallskip
\textit{Step~2: Well-posedness via the two-driver theory.}
Since $0\le y\ln(1+y)\le (\ln 2)\,y+y|\ln y|$ for $y>0$ and $0\le f_t^p\le\frac12|\theta_t|^2$, the driver $g_1^p$ satisfies
\[
0\le g_1^p(t,y,z)\le \bigl(\ln 2+\tfrac12|\theta_t|^2\bigr)y+y|\ln y|+\tfrac12\overline\lambda^2\,y|z|^2,
\]
hence Assumption~\hyperref[item:A']{A')} holds.
The map $z\mapsto g_2(t,y,z)=y\Lambda_t z$ is continuous, injective, with $|g_2(t,y,z)|\ge\underline\lambda\,y|z|$, so Assumption~(G2) holds.
The reduced driver
$\bar g^p(t,y,\gamma)=y\ln(1+y)+f_t^p\,y+\frac{|\gamma|^2}{2y}$
is jointly convex in $(y,\gamma)$ --- since $y\ln(1+y)$ is convex on $\R_+$, $f_t^p\,y$ is affine, and $|\gamma|^2/y$ is the perspective of $|\gamma|^2$ --- so Assumption~\hyperref[item:C']{C')} is satisfied.
By Proposition~\ref{prop:2DU}, the two-driver BSDE~\eqref{eq:sec7_twodriver_new} admits a unique solution $(Y^p,Z^p)$ with $0<Y^p\le (Y^p)^{h_1,h_2}$, $Y^p\in\mathcal{H}_T^r$ for all $r\ge1$, and $Z^p\in\mathcal{L}_T^2$.

Defining $\widetilde Z_t^p:=p_t+\Lambda_t Z_t^p$ and $V_t^p:=\ln(X_t^pY_t^p)$, It\^o's formula applied to $\ln(Y^p)$ together with~\eqref{eq:sec7_logwealth_new} shows that $(V^p,\widetilde Z^p)$ solves~\eqref{eq:sec7_primitive_bsde_new}. 
Uniqueness of $(V^p,\widetilde Z^p)$ follows from that of $(Y^p,Z^p)$: any other solution $(\widehat V^p,\widehat{\widetilde Z}^p)$ induces via $\widehat Y_t^p:=e^{\widehat V_t^p}/X_t^p$ a second solution to~\eqref{eq:sec7_twodriver_new}, which must coincide with $(Y^p,Z^p)$ by uniqueness. This proves item~(1).

\smallskip
\textit{Step~3: Optimality of $p^*=\theta$.}
For $p^1,p^2\in\A$ with $f_t^{p^1}\le f_t^{p^2}$ $d\mathbb{P}\times dt$-a.s., the reduced drivers satisfy $\bar g^{p^1}\le\bar g^{p^2}$, so Proposition~\ref{prop:2Dcomparisongen} yields $Y^{p^1}\le Y^{p^2}$ $d\mathbb{P}$-a.s.
Since $f_t^p\le f_t^\theta$ for all $p\in\A$, we obtain $Y^p\le Y^\theta$.
Evaluating at $t=0$ with $X_0^p=x$:
\[
V_0^p(X_T^p\xi)=\ln(xY_0^p)\le\ln(xY_0^\theta)=V_0^\theta(X_T^\theta\xi),
\]
hence $p^*=\theta$ is optimal and $\nu(x,\xi)=\ln(xY_0^\theta)$.
\end{proof}

\begin{remark}[Relation with \cite{BT19}]
The optimal control $p^\star=\theta$ coincides with the optimizer in \cite[Proposition~5.2]{BT19}. 
However, the opportunity equation there exhibits only the singular term $|z|^2/y$, whereas here it is a genuine two-driver BSDE containing both the logarithmic nonlinearity $y\ln(1+y)$ and the singular term.
Moreover, \cite{BT19} relies on a martingale optimality principle, while here optimality follows from the comparison theorem for two-driver BSDEs (Proposition~\ref{prop:2Dcomparisongen}).
Finally, we allow for unbounded terminal conditions $\xi$ under $L^p$-type assumptions, whereas \cite{BT19} requires $\xi$ bounded and bounded away from zero.
\end{remark}
\begin{remark}[Indifference pricing]
In the present setting, it is natural to define the indifference price of a contingent claim \(\xi\) in relative terms. More precisely, we define the \emph{indifference price factor} \(x_\xi>0\) by
\[
\nu\!\left(\frac{x}{x_\xi},\xi\right)=\nu(x,1), \qquad x>0,
\]
where \(1\) denotes the benchmark claim. By the previous analysis,
\[
\nu(x,\xi)=\ln x+\ln Y_0^\theta(\xi),
\]
hence
\[
x_\xi=\frac{Y_0^\theta(\xi)}{Y_0^\theta(1)}.
\]
As a consequence, the indifference price is a dimensionless multiplicative factor, independent of the initial endowment \(x\). Moreover, if \(\xi\ge 1\) then Proposition~\ref{prop:2Dcomparisongen} yields
\(
Y_0^\theta(\xi)\ge Y_0^\theta(1),
\)
and therefore \(x_\xi\ge 1\), consistently with the interpretation of \(x_\xi\) as the relative cost of the claim \(\xi\).
\end{remark}

\setcounter{equation}{0}\section{Applications to Return and Star-Shaped Risk Measures}\label{sec:appRRM}

In this section, we apply the results of the previous sections to study dynamic return and star-shaped risk measures
defined on $L^{p(2\delta+1)(e^B+1)}_+$. 
For further details on the terminology and concepts used in this section, we refer to Appendix~\ref{sec:drdm}.

We consider the GBSDE
\begin{equation}
\tilde\rho_t=X+\int_t^T\tilde\rho_s \tilde f(s,\tilde\rho_s,\tilde Z_s)\,ds
-\int_t^T\tilde \rho_s\tilde Z_s\,dW_s,
\label{eq:genGBSDEs}
\end{equation}
where $X\in L_+^{p(2\delta+1)(e^B+1)}$ and $\tilde f:\Omega\times[0,T]\times\R_+\times\R^n\to\R_+$ satisfies
\[
\tilde f(t,y,z)\leq \alpha_t+\beta_t|\ln(y)|+\delta |z|^2,\qquad d\mathbb{P}\times dt\text{-a.s.},\quad \forall(y,z)\in\R_+\times\R^n.
\]

\subsection{Structural properties}

The following theorem collects the main well-posedness and structural properties of the
dynamic (return) risk measures $(\tilde\rho_t)_{t\in[0,T]}$ induced by \eqref{eq:genGBSDEs}.
\begin{proposition}
\label{prop:return_gbsde}
Let $\alpha,\beta,\delta,p$ and $\tilde f$ satisfy the assumptions of
Corollary~\ref{cor:euGBSDE}. Then Equation~\eqref{eq:genGBSDEs} admits a unique positive
solution $(\tilde\rho,\tilde Z)$ with $\tilde\rho\in\mathcal{H}^{p(2\delta+1)(e^B+1)}_T$ and
$\tilde Z\in\mathcal{L}^2_T$. Moreover, for every $t\in[0,T]$, the mapping
$\tilde\rho_t:L^{p(2\delta+1)(e^B+1)}_+(\mathcal{F}_T)\to L^{p(2\delta+1)(e^B+1)}_+(\mathcal{F}_t)$
satisfies:
\begin{enumerate}
    \item \emph{(Monotonicity)} $\tilde\rho_t$ is monotone.
    \item \emph{(Lebesgue property)} If $(X^n)\subset L^{p(2\delta+1)(e^B+1)}_+$,
    $\sup_n X^n\in L^{p(2\delta+1)(e^B+1)}_+$, and $X^n\to X$ $\mathbb{P}$-a.s.,
    then $\tilde\rho_t(X^n)\to\tilde\rho_t(X)$ $\mathbb{P}$-a.s.
    \item \emph{(GA-convexity)} If $\tilde f$ satisfies
    \[
    \tilde f(t,y_1^{\lambda}y_2^{1-\lambda},\lambda z_1+(1-\lambda)z_2)
    \leq \lambda \tilde f(t,y_1,z_1)+(1-\lambda)\tilde f(t,y_2,z_2),
    \]
    then $\tilde\rho_t(X^{\lambda}Y^{1-\lambda})\leq\tilde\rho_t^{\lambda}(X)\tilde\rho_t^{1-\lambda}(Y)$
    for all $X,Y\in L_+^{p(2\delta+1)(e^B+1)}$ and $\lambda\in[0,1]$.
    \item \emph{(Positive homogeneity / Star-shapedness)} $\tilde\rho_t$ is positively
    homogeneous if $\tilde f$ does not depend on $y$, and star-shaped if $\tilde f$ is
    increasing in $y$.
    \item \emph{(Time-consistency)} The family $(\tilde\rho_t)_{t\in[0,T]}$ is time-consistent.
    \item \emph{(Normalization)} If $\tilde f(\cdot,1,0)\equiv 0$, then $\tilde\rho_t(1)=1$
    for all $t\in[0,T]$.
\end{enumerate}
\end{proposition}
\begin{remark}
GA-convexity in $(y,z)$ is strictly weaker than joint convexity in $(y,z)$
when $\tilde f$ is increasing in $y$ (see Example~\ref{ex:GAconv} below).
Positive homogeneity reflects scale-invariance in frictionless settings; star-shapedness
captures decreasing marginal impact of leverage and was introduced in \cite{CCMTW22,LR22}
(see also \cite{LRZ24}) and extended to the BSDE setting in \cite{LRZ23}.
Multiplicative convexity was first studied in \cite{BLR18,LR22} in the static case; here
we extend it to the dynamic setting. Time-consistency is naturally inherited from the flow
property of BSDEs (\cite{BEK05,DE92,CE02}).
\end{remark}

\subsection{Examples}

\begin{example}[Star-shaped risk measure]\label{ex:GAconv}
Let $\tilde{f}(t,y)=\beta_t\ln(1+y)$ with $\beta$ positive and bounded.
This driver is increasing in $y$, GA-convex (despite being concave in $y$),
and satisfies the required growth condition.
For $X\in L_+^{p(e^B+1)}$, $p>2$, $B=2T\|\beta\|_\infty^T$, the GBSDE
\[
\tilde\rho_t = X + \int_t^T \beta_s \tilde\rho_s \ln(1+\tilde\rho_s)\,ds
- \int_t^T \tilde\rho_s \tilde Z_s\,dW_s
\]
admits a unique solution, and $(\tilde\rho_t)_{t\in[0,T]}$ is monotone, star-shaped,
multiplicatively convex, and time-consistent by Proposition~\ref{prop:return_gbsde}.
By analogy with the monetary case, the parameter $\beta$
can be interpreted as a stochastic discount rate encoding ambiguity on interest rates (see \cite{ELKR09,LRZ24}),
with star-shapedness replacing positive homogeneity in the geometric setting.
\end{example}

\begin{example}[$\gamma$-norms and robust $\gamma$-norms]\label{ex:gamma_norms}
Fix $\gamma>1$ and $X\in L^{p\gamma}_+$.

\medskip
\noindent\emph{\textbf{$\gamma$-norms.}} The GBSDE
\begin{equation}
\tilde\rho_t=X+\int_t^T\frac{\gamma-1}{2}\tilde\rho_s|\tilde Z_s|^2\,ds
-\int_t^T\tilde\rho_s\tilde Z_s\,dW_s
\label{eq:gammanorms}
\end{equation}
has a unique positive solution $(\tilde\rho,\tilde Z)\in\mathcal{H}^{p\gamma}_T\times\mathcal{L}^2_T$.
Applying It\^o's formula to $\tilde\rho_{\fatdot{}}^{\gamma}$ yields
$\tilde\rho_t^{\gamma}=\mathbb{E}[X^{\gamma}|\mathcal{F}_t]$,
so $\tilde\rho_t=\mathbb{E}[X^{\gamma}|\mathcal{F}_t]^{1/\gamma}$ is the \emph{dynamic $\gamma$-norm} (cf.\ \cite{BLR18,BLR21,W08}),
and uniqueness follows from that of the resulting martingale equation.

\medskip
\noindent\emph{\textbf{Robust $\gamma$-norms.}}
Let $g:[0,T]\times\Omega\times\R^n\to\R_+$ be $\mathcal{P}\times\mathcal{B}(\R^n)$-measurable,
convex, positively homogeneous, and satisfy $0\leq g(t,z)\leq C|z|$ a.s.\ for some $C>0$.
For $p\geq2$, the GBSDE
\begin{equation}
\tilde\rho_t=X+\int_t^T\tilde\rho_s\!\left(g(s,\tilde Z_s)+\frac{\gamma-1}{2}|\tilde Z_s|^2\right)ds
-\int_t^T\tilde\rho_s\tilde Z_s\,dW_s
\label{eq:robustgammanorms}
\end{equation}
has a unique positive solution. Setting $\tilde V:=\gamma\tilde\rho^{\gamma}\tilde Z$ and
applying It\^o's formula to $\tilde\rho_{\fatdot{}}^{\gamma}$ gives a BSDE with sublinear
driver $g$, whose solution satisfies (see \cite{BEK05,DKRT14})
\[
\tilde\rho_t=\sup_{\mu\in\mathcal{A}}\mathbb{E}_{\mathbb{Q}^{\mu}}\!\left[X^{\gamma}\big|\mathcal{F}_t\right]^{1/\gamma},
\qquad d\mathbb{P}\times dt\text{-a.s.},
\]
where $\mathcal{A}:=\{(\mu_t): |\mu_t|\leq C\text{ a.s.}\}$ and $\mathbb{Q}^{\mu}$
has density $\mathcal{E}^{\mu}_T$. Thus $\tilde\rho_t$ is the \emph{dynamic robust $\gamma$-norm},
the natural return counterpart of entropy coherent risk measures (\cite{CF10,LS13}).

\medskip
In both cases, their unique solution, monotonicity, positive homogeneity, multiplicative convexity, and
time-consistency follow directly from Proposition~\ref{prop:return_gbsde}.
\end{example}

\noindent\textbf{Acknowledgements.} {\small We are 
grateful to 
Fabio Bellini, Sonja Cox, Freddy Delbaen, Marco Frittelli, Michael Kupper, Antonis Papapantoleon, Mitja Stadje and to conference participants at the 2024 Bachelier Finance Society World Congress in Rio de Janeiro, the 2024 Probability Conference in Rome and the 2025 
General AMaMeF Conference in Verona for useful comments and discussions. 
This research was funded in part by the Netherlands Organization for Scientific Research under grant NWO Vici 2020--2027 (Laeven) and by an Ermenegildo Zegna Founder's Scholarship (Zullino). 
Rosazza Gianin and Zullino are members of Gruppo Nazionale per l’Analisi Matematica, la Probabilità e le loro Applicazioni (GNAMPA, INdAM), Italy, and acknowledge the financial support of Gnampa Research Project 2024 (PRR-20231026-073916-203).}

\newpage

{\Huge
\begin{center}
SUPPLEMENTARY MATERIAL TO\\
``Geometric BSDEs''\\
(FOR ONLINE PUBLICATION)
\end{center}
}

\newpage

\appendix
\renewcommand{\thesection}{\Alph{section}}
\section{Online Appendix}
 \label{app}
 \renewcommand{\thesubsection}{\thesection.\Roman{subsection}}
\renewcommand{\theequation}{\thesection.\Roman{equation}}
\renewcommand{\thetheorem}{\thesubsection.\arabic{theorem}}
\renewcommand{\thetheorem}{\thesection.\Roman{theorem}}
\setcounter{equation}{0}
\renewcommand{\thesection}{\Alph{section}}
\renewcommand{\thesubsection}{\thesection.\arabic{subsection}}
\renewcommand{\thetheorem}{\thesubsection.\arabic{theorem}}
\setcounter{theorem}{0}

\subsection{Additional material to Section~\ref{sec:prel}}
\begin{lemma}
    For any $t\in[0,T]$, the space $\mathcal{L}^{\infty}(\mathcal{F}_t)$ can be identified with the set of random variables given by
    \begin{equation*}L^{\infty}_{\ln}(\mathcal{F}_t):=\{Y\in L^{0}_{+}: \ \ln(Y)\in L^{\infty}(\mathcal{F}_t)\}.\end{equation*}\label{lem:spaceidentification}\end{lemma}
\begin{proof}
\vspace{-0.4cm}The proof is a routine verification.
\end{proof}
\subsubsection{Dynamic return and monetary risk measures} 
\label{sec:drdm}
In this subsection, we provide definitions of dynamic return and monetary risk measures, used throughout. 
\begin{definition}\label{def:drm}
    Let $L^{\infty}(\mathcal{F}_T)\subseteq\mathcal{X}\subseteq L^{0}(\mathcal{F}_T)$ and let $t \in [0,T]$. 
    Then $\rho_t:\mathcal{X}(\mathcal{F}_T)\to \mathcal{X}(\F_t)$ is a dynamic risk measure if it is monotone, i.e., for any $X,Y\in \mathcal{X}(\mathcal{F}_T)$ such that $X\geq Y$, it holds that $\rho_t(X)\geq\rho_t(Y)$.
    A dynamic monetary risk measure is a dynamic risk measure that additionally verifies cash-additivity, that is, for any $t \in [0,T]$, $X \in \mathcal{X}(\F_T)$ and $\eta_t \in \mathcal{X}(\F_t)$ such that $X+\eta_t\in\mathcal{X}(\mathcal{F}_T)$, it holds that $\rho_t( X + \eta_t)= \rho_t(X)+ \eta_t$.
    
    Let $\mathcal{L}^{\infty}(\mathcal{F}_T)\subseteq\mathcal{Y}\subseteq L^0_{+}(\mathcal{F}_T)$. 
    Then $\tilde\rho_t:\mathcal{Y}(\F_T)\to \mathcal{Y}(\F_t)$ is a dynamic return risk measure if it is monotone on $\mathcal{Y}(\mathcal{F}_T)$ and positively homogeneous with respect to $\F_t$-measurable random variables, that is, for any $t \in [0,T]$, $X \in \mathcal{X}(\F_T)$ and $\xi_t \in \mathcal{L}^{\infty}(\F_t)$ such that $\xi_t\cdot X\in\mathcal{Y}(\mathcal{F}_T)$, it holds that $\tilde{\rho}_t( \xi_t \cdot X)= \xi_t \cdot \tilde{\rho}_t(X)$.
   \end{definition}
By a slight abuse of notation, we will use the term \emph{dynamic risk measure} to refer interchangeably to a single map \(\rho_t\) (resp.\ $\tilde\rho_t)$ and to the family \((\rho_t)_{t\in[0,T]}\) (resp.\ $(\tilde\rho_t)_{t\in[0,T]}$), the intended meaning being clear from the context.

   We present a (non-exhaustive) list of axioms for dynamic monetary and return risk measures; see, e.g.,  \cite{D12,FS11,BLR18,LR22} for further details and discussion. 
   Let $t\in[0,T]$.
\begin{definition}\label{def:prop}
A dynamic risk measure $\rho_t:\mathcal{X}(\mathcal{F}_T)\to\mathcal{X}(\mathcal{F}_t)$ is: (i) convex if for any $X,Y\in\mathcal{X}(\mathcal{F}_T)$ and $\lambda\in[0,1]$ such that $\lambda X+(1-\lambda)Y\in\mathcal{X}(\mathcal{F}_T)$ it holds that $\rho_t(\lambda X+(1-\lambda)Y)\leq\lambda\rho_t(X)+(1-\lambda)\rho_t(Y)$; (ii) positively homogeneous if for any $X\in\mathcal{X}(\mathcal{F}_T)$ and $\eta_t\in L^{\infty}_+(\mathcal{F}_t)$ such that $\eta_t\cdot X\in\mathcal{X}(\mathcal{F}_T)$ it results that $\rho_t(\eta_t\cdot X)=\eta_t\cdot\rho_t(X).$
    
A monotone functional $\tilde\rho_t:\mathcal{Y}(\mathcal{F}_T)\to\mathcal{Y}(\mathcal{F}_t)$ is: (i) multiplicatively convex if for any $X,Y\in\mathcal{Y}(\mathcal{F}_T)$ and $\lambda\in[0,1]$ such that $X^{\lambda}Y^{1-\lambda}\in\mathcal{Y}(\mathcal{F}_T)$ it holds that $\rho_t(X^{\lambda}Y^{1-\lambda})\leq\rho^{\lambda}_t(X)\rho^{1-\lambda}_t(Y)$; (ii) star-shaped if for any $X\in\mathcal{Y}(\mathcal{F}_T)$ and $\eta_t\in {\mathcal{L}}^{\infty}_+(\mathcal{F}_t)$ with $\eta_t\leq 1$ $d\mathbb{P}$-a.s.\ such that $\eta_t\cdot X\in\mathcal{Y}(\mathcal{F}_T)$ it results that $\rho_t(\eta_t\cdot X)\leq\eta_t\cdot\rho_t(X).$

A (family of) risk measure(s) $\rho_t:\mathcal{X}(\mathcal{F}_T)\to\mathcal{X}(\mathcal{F}_t)$ for $t\in[0,T]$ is time-consistent on $\mathcal{X}(\mathcal{F}_T)$ if $\rho_s(X)=\rho^{t}_s\big(\rho_t(X)\big)$ for any $t,s\in[0,T]$ with $s<t$ and $X\in\mathcal{X}(\mathcal{F}_T)$ where $\rho^t_s: \mathcal{X}(\mathcal{F}_t)\to\mathcal{X}(\mathcal{F}_s)$ is the restriction of $\rho_s$ to $\mathcal{X}(\mathcal{F}_t)$, i.e., $\rho_s^t := (\rho_s)_{|\mathcal{X}(\mathcal{F}_t)}$. 
Time-consistency of functionals defined on $\mathcal{Y}(\mathcal{F}_T)$, such as $\tilde\rho_t$, can be formulated analogously.
\end{definition}

We presented the preceding definition without imposing a linear structure on $\mathcal{X}$ and $\mathcal{Y}$, recognizing that certain spaces we will utilize hereafter may lack this property. 
A one-to-one correspondence between return risk measures and monetary risk measures has been proved in \cite{BLR18} in a static setting, using as reference spaces $\mathcal{X}\equiv L^{\infty}$ and $\mathcal{Y}\equiv\mathcal{L}^{\infty}$. 
We can extend this one-to-one correspondence to the dynamic case.
Starting from a dynamic monetary risk measure $\rho_t:L^{\infty}(\mathcal{F}_T) \to L^{\infty}(\mathcal{F}_t)$, there is a natural way to define the corresponding dynamic return risk measure
{$\tilde\rho_t: \mathcal{L}^{\infty}(\mathcal{F}_T)\to \mathcal{L}^{\infty}(\mathcal{F}_t)$} by setting
\begin{equation}
    \tilde\rho_t(X):=\exp(\rho_t(\ln(X)).
    \label{eq: rho tilde from rho}
\end{equation}
We notice that, for each $t\in[0,T]$, $\tilde\rho_t$ is well-defined by Lemma~\ref{lem:spaceidentification}. 
Indeed, if $X\in \mathcal{L}^{\infty}(\mathcal{F}_T)$, then $\ln(X)\in L^{\infty}(\mathcal{F}_T)$, thus $\rho_t(\ln(X))\in {L}^{\infty}(\mathcal{F}_t)$, hence $\exp(\rho_t(\ln(X))\in \mathcal{L}^{\infty}(\mathcal{F}_t).$
Conversely, given a dynamic return risk measure $\tilde\rho_t: \mathcal{L}^{\infty}(\mathcal{F}_T)\to \mathcal{L}^{\infty}(\mathcal{F}_t)$, we can define the corresponding dynamic monetary risk measure $\rho_t:{L}^{\infty}(\mathcal{F}_T) \to L^{\infty}(\mathcal{F}_t)$ via the formula
\begin{equation}
\rho_t(X):=\ln \left(\tilde\rho_t \left(e^X \right)\right).
\label{eq: rho from rho tilde}
\end{equation}
Once again, $\rho_t$ is well defined as $e^{X}\in\mathcal{L}^{\infty}(\mathcal{F}_T)$ for any $X\in L^{\infty}(\mathcal{F}_T)$, thus $\tilde\rho_t(e^X)\in \mathcal{L}^{\infty}(\mathcal{F}_t)$ and Lemma~\ref{lem:spaceidentification} ensures that $\ln(\tilde\rho_t(X))\in L^{\infty}(\mathcal{F}_t).$ 
Properties of $\tilde\rho_t$ and $\rho_t$ are summarized in the following proposition.

\begin{proposition} For any $t\in[0,T]$, let $\tilde\rho_t$ and ${\rho}_t$ be as defined in \eqref{eq: rho tilde from rho} and \eqref{eq: rho from rho tilde}. Then:

\noindent (a) $\rho_t(0)=0$ (normalization for $\rho$) if and only if $\tilde{\rho}_t(1)=1$ (normalization for $\tilde{\rho}$) for any $t\in[0,T].$

\noindent (b) $\rho_t$ satisfies cash-additivity on $L^{\infty}(\mathcal{F}_T)$ if and only if $\tilde{\rho}_t$ satisfies positive homogeneity on $\mathcal{L}^{\infty}(\F_T)$.

\noindent (c) $\rho_t$ satisfies cash-super-additivity on $L^{\infty}(\mathcal{F}_T)$, i.e., $\rho_t(X+m_t)\geq\rho_t(X)+m_t$ for any $X\in L^{\infty}(\mathcal{F}_T)$ and $m_t\in L_+^{\infty}(\mathcal{F}_t)$, if and only if $\tilde{\rho}_t$ satisfies star-shapedness on $\mathcal{L}^{\infty}(\F_T)$.

\noindent (d) $\rho_t$ is monotone on $L^{\infty}(\mathcal{F}_T)$ if and only if $\tilde{\rho}_t$ is monotone on $\mathcal{L}^{\infty}(\F_T)$.

\noindent (e) $\rho_t$ is convex on  $L^{\infty}(\mathcal{F}_T)$ if and only if $\tilde{\rho}_t$ is multiplicatively convex on $\mathcal{L}^{\infty}(\F_T)$.

\noindent (f) $\rho_t$ is positively homogeneous on $L^{\infty}(\mathcal{F}_T)$ if and only if $\tilde \rho_t$ is multiplicatively positively homogeneous on $\mathcal{L}^{\infty}(\F_T)$, i.e., for any $t \in [0,T]$, $X \in \mathcal{L}^{\infty}(\F_T)$ and $\eta_t \in \mathcal{L}^{\infty}(\mathcal{F}_t)$, $\tilde{\rho}_t( X^{\eta_t} ) = \tilde{\rho}_t^{\eta_t} (X)$.

\noindent (g) The family of functionals $(\rho_t)_{t\in[0,T]}$ is time-consistent on $L^{\infty}(\mathcal{F}_T)$ if and only if $(\tilde{\rho}_t)_{t\in[0,T]}$ is time-consistent on $\mathcal{L}^{\infty}(\F_T)$.
\label{prop:1to1P}
\end{proposition}
\begin{proof}
Items (a)-(g) can easily be obtained by applying \eqref{eq: rho tilde from rho} and \eqref{eq: rho from rho tilde}.
\end{proof}

\subsection{Proofs} \label{app: proofs}
\setcounter{theorem}{0}

\begin{proof}[Proof of Proposition~\ref{prop:equivalence}]
We start by proving the first statement. 
Fix $X\in L^{0}_{+}(\mathcal{F}_T)$ and define $f(x):=\ln(x)$. 
Applying It\^o's formula to $ Y'_{\fatdot{}}=f(Y_{\fatdot{}})$ we obtain
\begin{align*}
-d Y'_t=e^{- Y'_t}\left(g(t,e^{Y'_t},Z_t)+\frac{e^{-Y'_t}|Z_t|^2}{2}\right)-e^{- Y'_t}Z_tdW_t.
\end{align*}
Setting $Z'_t:=e^{-Y'_t}Z_t$ it holds that
\begin{align*}
-dY'_t=e^{-Y'_t}\left(g(t,e^{Y'_t},e^{Y'_t}Z'_t)+\frac{e^{Y'_t}|Z'_t|^2}{2}\right)-Z'_tdW_t,
\end{align*}
with terminal condition $X'=\ln(X)$. 
By defining $g'(t,y,z):=e^{-y}g(t,e^y,e^yz)+\frac{|z|^2}{2}$, it results that
$$|g'(t,y,z)|\leq e^{-y}\left(\alpha_t e^y+\beta_t e^y|\ln(e^y)|+\delta\frac{e^{2y}|z|^2}{e^y}\right)+\frac{|z|^2}{2}=\alpha_t+\beta_t|y|+\Big(\delta+\frac{1}{2}\Big)|z|^2,$$
thus the BSDE with parameters $(X',g')$ admits a solution $(Y', Z'):=(\ln(Y),Z/Y)$, in the sense of Definition~\ref{def:solBSDE}. 
Indeed, $Y$ is a positive and continuous process, as composition of continuous functions in their domain, $\int_0^T|Z'_s|^2ds\leq\frac{1}{\essinf_t Y_t}\int_0^T|Z_s|^2ds<+\infty$ and $$\int_0^T|g'(s,Y'_s,Z'_s)|ds\leq \int_0^T\alpha_s+\beta_s|Y'_s|+\Big(\delta+\frac{1}{2}\Big)|Z'_s|^2ds<+\infty.$$ 
 
The other implication can be proved analogously by employing the substitution $Y:=\exp(Y')$, verifying that $(Y,Z)=(\exp(Y'),\exp(Y')Z')$ is a solution to the BSDE \eqref{eq:LNQ} with parameters $(X,g)$, where $X:=\exp(X')$.
\end{proof}

\begin{proof}[Proof of Corollary~\ref{cor:gen1to1}]
Since $X\in\mathcal{L}^{\infty}(\mathcal{F}_T)$, $X':=\ln(X)\in L^{\infty}(\mathcal{F}_T)$, hence the corresponding quadratic BSDE with parameters $(X',g')$ with $g'$ as in Proposition~\ref{prop:equivalence} admits a maximal and minimal solution verifying $(Y',Z')\in\mathcal{H}^{\infty}_T\times\text{BMO}(\mathbb{P})$, by results in \cite{K00}. 
Thus, $(Y,Z):=(\exp(Y'),\exp(Y')Z')\in\mathcal{H}^{\infty}_T\times\text{BMO}(\mathbb{P})$ is a solution to the BSDE with parameters $(X,g)$.    
\end{proof}

\begin{proof}[Proof of Proposition~\ref{prop:EUgen}]
\textit{(i)} Consider
\begin{align}
    -d\bar Y_t
    =\left((2\delta+1)\alpha_t\bar Y_t+\beta_t \bar Y_t|\ln(\bar Y_t)|\right)dt-\bar Z_tdW_t,
    \label{eq:posol}
\end{align}
with terminal condition $\bar X= X^{2\delta+1}$. 
We define $$\bar h(t,y,z):=(2\delta+1)\alpha_t y+\beta_t y|\ln(y)|, \ \ d\mathbb{P}\times dt\text{-a.s., } \forall (y,z)\in\R_+\times\R^n.$$ 
Equation~\eqref{eq:posol} admits a positive solution if and only if Equation~\eqref{eq:LNQ} does.
To see this, it is enough to apply It\^o's formula to $\bar Y=u(Y)$ to get Equation~\eqref{eq:posol} from Equation~\eqref{eq:LNQ} and \textit{vice versa} for $Y=u^{-1}(\bar Y)$, where $u(x):=x^{2\delta+1}$ for any $x>0$. 

The existence of a positive solution to Equation~\eqref{eq:posol} is equivalent to proving the existence of a positive solution $(w,\zeta)$ to
\begin{equation}
w_t=\ln(1+\bar X)+\int_t^T\tilde h(s,w_s,\zeta_s)ds-\int_t^T\zeta_sdW_s,
\label{eq:2posol}
\end{equation}
where, $d\mathbb{P}\times dt\text{-a.s., }  \forall(y,z)\in\R_+\times\R^n$, $$\tilde h(t,y,z):=((2 \delta +1)\alpha_s(e^y-1)+\beta_s(e^y-1)|\ln(e^y-1)|)e^{-y}+\frac{1}{2}|z|^2.$$ 
To see this, it suffices to apply It\^o's formula to $w:=\ln(1+\bar Y)$ and $\bar Y = \exp(w)-1.$
Since $\frac{x|\ln(x)|}{x+1}\leq 1+|\ln(x+1)|$, it holds that
$$0\leq\tilde h(t,y,z)\leq (2 \delta +1) \alpha_t+\beta_t+\beta_ty+\frac{1}{2}|z|^2, \ \ d\mathbb{P}\times dt\text{-a.s., }  \forall(y,z)\in\R_+\times\R^n.$$
As $\tilde X:=\ln(1+\bar X)>0$, by Theorem~3.1~(i) of \cite{B20}, there exists a positive solution to the Equation~\eqref{eq:2posol} with parameters $(\tilde X,\tilde h)$, yielding the thesis. 
In addition, the regularity of $w$ provided by Theorem~3.1~(i) yields the required regularity for $Y=(\exp(w)-1)^{\frac{1}{2\delta+1}}$. 
\smallskip

\noindent\textit{(ii)} Now we assume (H1)+(H2)'. 
Clearly, existence of the solution follows as above, while the regularity of $Y$ is a consequence of Proposition~3.2 in \cite{B20}. 
We prove that $Z\in\mathcal{M}^2_T$. We can use a similar technique as in Proposition~3.5 of \cite{BT19}.

Let us start by considering the case $\delta\neq\frac{1}{2}$. 
We apply It\^o's formula to $f(x)=x^2$, obtaining
\begin{align*}
    Y_t^2=X^2+\int_t^T \left(2\alpha_s Y_s^2+2\beta_s Y_s^2|\ln(Y_s)|+(2\delta-1)|Z_s^2| \right)ds-2\int_t^TY_sZ_sdW_s.
\end{align*}
Considering the sequence of stopping times $\tau_n:=\inf\{t>0:\int_0^t4|Y_s|^2|Z_s|^2ds\geq n\}$, rearranging and taking the expectation, it holds that
\begin{align*}
    &|2\delta-1|\mathbb{E}\left[\int_0^{T\wedge\tau_n}|Z_t|^2dt\right]\leq \mathbb{E}\left[X^2+Y_0^2+\int_0^{T\wedge \tau_n} \left(2\alpha_t Y_t^2+2\beta_t Y_t^2|\ln(Y_t)| \right) dt\right] \\
    &\leq\mathbb{E}\left[X^2+Y_0^2+\int_0^{T\wedge \tau_n} \left(\frac{2\alpha_t^q}{q}+\frac{2Y_t^{2p}}{p}+\frac{2\beta_t^q}{q}+\frac{2Y_t^{2p}|\ln(Y_t)|^p}{p} \right)dt\right],
\end{align*}
where the last inequality follows from Young's inequality.
In order to study the integrability of $Z$, we make use of the inequality $y^{2p}|\ln(y)|^p\leq K_{\varepsilon}+y^{2p+\varepsilon}$ that holds for any $\varepsilon>0$, with $K_{\varepsilon}>0$ depending only on the chosen $\varepsilon>0$. This yields
$$|2\delta-1|\mathbb{E}\left[\int_0^{T\wedge\tau_n}|Z_t|^2dt\right]\leq K_{\varepsilon}\mathbb{E}\left[\int_0^{T} \left(\frac{2\alpha_t^q}{q}+\frac{2\beta_t^q}{q} \right)dt+K_T\sup_{t\in[0,T]}\left\{Y_t^{2p+\varepsilon}\right\}\right],$$ where $K_T>0$. 
Letting $n\to\infty$, we have $T\wedge \tau_n\to T$ $d\mathbb{P}$-a.s., and Fatou's lemma yields
$$|2\delta-1|\mathbb{E}\left[\int_0^{T}|Z_t|^2dt\right]\leq K_{\varepsilon} \mathbb{E}\left[\int_0^{T} \left(\frac{2\alpha_t^q}{q}+\frac{2\beta_t^q}{q} \right)dt+K_T\sup_{t\in[0,T]}\left\{Y_t^{2p+\varepsilon}\right\}\right].$$ 
Choosing $\varepsilon=2\delta p>0$, we obtain the thesis, since $2p+\varepsilon=2p(\delta+1)\leq p(2\delta+1)(e^{\inf_{\omega\in\Omega}B}+1)$. 

If $\delta=\frac{1}{2}$, It\^o's formula applied to the function $x^2\ln(x)$ entails
\begin{align*}
    Y_t^2\ln(Y_t)=\:&X^2\ln(X)\\
    &+\int_t^T \left(\alpha_s Y_s^2(1+2 \ln(Y_s))+\beta_s Y_s^2|\ln(Y_s)|(1+ 2\ln(Y_s))-|Z_t|^2 \right)ds \\
    &-\int_t^T(2Y_s\ln(Y_s)+Y_s)Z_sdW_s.
\end{align*}
Considering the sequence of stopping times 
$\tau_n:=\inf\{t>0:\int_0^t|Y_sZ_s|^2+|Y_s\ln(Y_s)Z_s|^2ds\geq n\}$ and arguing similarly as above we get
$$\mathbb{E}\left[\int_0^T|Z_t|^2dt\right]\leq K\mathbb{E}\left[\int_0^{T} \left(\frac{2\alpha_t^q}{q}+\frac{2\beta_t^q}{q} \right)dt+\sup_{t\in[0,T]}\left\{Y_t^{2p(\delta+1)}\right\}\right]<+\infty,$$
for some $K>0.$\smallskip

\noindent\textit{(iii)} Finally, consider (H1)'+(H2)''. 
Following the same passages as above and including the new term $\gamma_t|z|$ in the driver, we obtain the transformed equation
\begin{equation}
w'_t=\ln(1+\bar X)+\int_t^T h'(s,w'_s,\zeta'_s)ds-\int_t^T\zeta'_sdW_s,
\label{eq:3posol}
\end{equation}
where
\begin{align*}
 h'(t,y,z)&:=((2 \delta +1)\alpha_s(e^y-1)+\beta_s(e^y-1)|\ln(e^y-1)|)e^{-y}+\gamma_t|z|+\frac{1}{2}|z|^2 \\
&\leq (2 \delta +1) \alpha_t+\beta_t+\beta_ty+\gamma_t|z|+\frac{1}{2}|z|^2.
\end{align*}
By the proof of Proposition~3.3 in \cite{B20}, Equation~\eqref{eq:3posol} admits a positive solution if the BSDEs with parameters $(\xi,h'')=(\exp((\bar X+\int_0^T\alpha_s+\beta_sds)(e^{\int_0^T{\beta_sds}}+1))-1,\gamma_{\fatdot{}}|z|)$ admits a positive solution $(Y',Z')$. 
By assumption (H2)'' and Theorem~2.1 in \cite{BEH15}, the BSDE with parameters $(\xi,h'')$ admits a (unique) solution $(Y',Z')\in\mathcal{H}^{p}_T\times\mathcal{M}^p_T$. 
In addition, $Y'>0$ since $\xi>0$. 
Indeed, after choosing a suitable localization $(\tau'_n)_{n\in\mathbb{N}}$ as above, it holds that
$$Y'_t=\mathbb{E}\left[\left.\xi+\int_t^{T\wedge\tau'_n}\gamma_s|Z'_s|ds\right|\mathcal{F}_t\right]\geq\mathbb{E}[\xi|\mathcal{F}_t]>0,$$ 
where we employed the positivity of $\gamma|Z'|.$
Thus, the thesis follows with the required regularity for the $y$-component of the solution. 

It remains to be shown that $Z\in\mathcal{M}^2_T.$ 
Suppose $\delta\neq1/2$ and apply It\^o's formula to $f(x)=x^2$, which yields
\begin{align*}
    Y_t^2=\ &X^2+\int_t^T \left(2\alpha_s Y_s^2+2\beta_s Y_s^2|\ln(Y_s)|+2\gamma_sY_s|Z_s|+(2\delta-1)|Z_s^2| \right)ds\\
    &-2\int_t^TY_sZ_sdW_s.
\end{align*}
Taking $\tau_n$ as above and performing similar calculations, it holds that
\begin{align*}
    &|2\delta-1|\mathbb{E}\left[\int_0^{T\wedge\tau_n}|Z_t|^2dt\right]\\
    &\leq \mathbb{E}\left[X^2+Y_0^2+\int_0^{T\wedge \tau_n} \left(2\alpha_t Y_t^2+2\gamma_tY_t|Z_t|+2\beta_t Y_t^2|\ln(Y_t)|\right)dt\right] \\
    &\leq\mathbb{E}\Bigg[X^2+Y_0^2\\
    &\qquad+\int_0^{T\wedge \tau_n} \left(\frac{2\alpha_t^q}{q}+\frac{2Y_t^{2p}}{p}+\frac{2\beta_t^q}{q}+\frac{2Y_t^{2p}}{p\lambda}+\frac{2\gamma_t^{2q}}{q\lambda}+\frac{\lambda Z_t^{2}}{2}+\frac{2Y_t^{2p}|\ln(Y_t)|^p}{p} \right)dt\Bigg],
\end{align*}
where we have used Young's inequality. 
Then, the thesis follows arguing as above by choosing $\lambda=|2\delta-1|.$ 
The case $\delta=1/2$ can be shown similarly, employing the substitution $f(x)=x^2\ln(x)$.
\end{proof}

We need a preliminary lemma, provided in \cite{B20}. 
For the reader's convenience, we state this result again.
\begin{lemma}
\label{lem:bahlali}
Let $g,g_1,g_2:\Omega\times[0,T]\times\R\times\R^n\to\R$ be $\mathcal{P}\times\mathcal{B}(\R)\times \mathcal{B}(\R^n)/\mathcal{B}(\R)$-measurable functions and $X,X_1,X_2\in L^0(\mathcal{F}_T)$ with $X_1\leq X\leq X_2$. 
Let us assume that the BSDEs with parameters $(X_1,g_1)$ and $(X_2,g_2)$ admit solutions $(Y^1,Z^1)$ and $(Y^2,Z^2)$, respectively, such that:
\begin{itemize}
    \item $Y^1\leq Y^2$,
    \item $d\mathbb{P}\times dt$-a.s., $y\in\left[Y_t^1(\omega),Y_t^2(\omega)\right]$ and $z\in\R^n$ it holds that:
    $g_1(t,y,z)\leq g(t,y,z)\leq g_2(t,y,z)$ and $|g(\omega,t,y,z)|\leq \eta_t(\omega)+C_t(\omega)|z|^2$, where $C$ and $\eta$ are pathwise continuous and $(\mathcal{F}_t)_{t\in[0,T]}$-adapted processes with $\eta$ verifying $\int_0^T|\eta_s|ds<+\infty$ $d\mathbb{P}$-a.s.
\end{itemize}
Then, if $g$ is continuous in $(y,z)$ $d\mathbb{P}\times dt$-a.s., the BSDE with parameters $(X,g)$ admits at least one solution $(Y,Z)$ in the sense of Definition~\ref{def:solBSDE}, with $Y^1\leq Y\leq Y^2.$ 
In addition, among all solutions lying between $Y^1$ and $Y^2$, there exist a maximal and a minimal solution.
\end{lemma}

\begin{proof}[Proof of Corollary~\ref{cor:EU1}]
    With the same notation as in Proposition~\ref{prop:EUgen}, we know that the BSDE with parameters $(\bar X,\bar h)$ admits a unique solution in $\mathcal{H}^{e^{\lambda T}+1}_T\times\mathcal{M}^2_T$, according to Theorem~2.1 in \cite{BKK17}. 
    Hence, $Y\in\mathcal{H}_T^{(2\delta+1)(e^{\lambda T}+1)}$, while the regularity of $Z$ can be checked as in the proof of Proposition~\ref{prop:EUgen}. 
    Thus, considering two solutions $(Y,Z)$ and $(U,Q)$ to the BSDE with parameter $(X,h)$, we have $\bar Y=Y^{2\delta+1}=U^{2\delta+1}=\bar U\in\mathcal{H}^{(2\delta+1)(e^{\lambda T}+1)}_T$ and $\bar Z=(2\delta+1)Y^{2\delta}Z=(2\delta+1)U^{2\delta}Q=\bar Q\in\mathcal{M}^2_T$. 
    In particular, $U=Q$ $dt\times d\mathbb{P}\mbox{-a.s.}$ since
\begin{align*}\int_0^T|Z_t-Q_t|^2dt&=\int_0^TY^{-4\delta}_t|Y^{2\delta}_tZ_t-U^{2\delta}_tQ_t|^2dt \\ &\leq \esssup_{t\in[0,T]}\{Y^{-4\delta}_t\}\int_0^T|Y^{2\delta}_tZ_t-U^{2\delta}_tQ_t|^2dt\\
&\leq C_{\omega}\int_0^T|\bar Z_t-\bar Q_t|^2dt=0 \ d\mathbb{P}\mbox{-a.s.},
\end{align*}
where the equality $Y=U$ leads to the above equality, while the last inequality is due to the pathwise continuity of $Y$, $Y>0$ and $\bar Z= \bar Q$.
\end{proof}

\begin{proof}[Proof of Proposition~\ref{prop:Linf}]
The existence and uniqueness of the solution to Equation~\eqref{eq:exdr} has already been established in Corollary~\ref{cor:EU1}. 
Analogously, Theorem~\ref{th:Egen} guarantees the existence of a minimal and a maximal solution to the BSDE with parameters $(X,g)$. 
We need to check the regularity of these solutions. \smallskip

\noindent\textit{Boundedness of $Y$ for Equation~\eqref{eq:exdr}:} 
It suffices to show that $\bar Y^h:=(Y^h)^{2\delta+1}$ is bounded, as this implies the boundedness of $Y^h$. 
We can prove that $\bar Y^h$ is bounded by employing the same transformation as in Proposition~\ref{prop:EUgen} and applying Corollary~3.4 (i) in \cite{B20}.
\smallskip

\noindent\textit{Boundedness of $Y$ for a general driver $g$:} 
By Theorem~\ref{th:Egen}, we know that the $Y$-component of the solution to the BSDE with parameters $(X,g)$ verifies $0<Y\leq Y^{h}$, hence boundedness of $Y^{h}$ implies the same property for $Y$. 

\noindent\textit{$Z\in\text{BMO}(\mathbb{P}):$} 
In order to prove that $Z\in\text{BMO}(\mathbb{P})$ we first show that $Z\in\mathcal{M}^2_T$. 
Proceeding as in the proof of Theorem~\ref{th:Egen}, applying It\^o's formula to $f(x):=\ln(1+x)$, it holds that
$$\frac{1}{2(1+\|Y\|^T_{\infty})^2}\mathbb{E}\left[\int_0^{T}|Z_s|^2ds\right]\leq\mathbb{E}\left[\int_0^{T}\frac{1}{2(1+Y_s)^2}|Z_s|^2ds\right]\leq Y_0+\|X\|_{\infty}<+\infty,$$
where we employed boundedness of $X$. 
Hence, $Z\in\mathcal{M}^2_T.$
Having proved that $Z\in\mathcal{M}^2_T$ and $Y$ is bounded and positive, the process $\frac{Z}{1+Y}\in\mathcal{M}^2_T$. 
Thus, we can follow the passage as above, without incorporating the localization procedure, as the stochastic integral of $\frac{1}{(1+Y)^2}|Z|^2$ is a uniformly integrable martingale. 
Consequently, we have
\begin{align*}
&\frac{1}{2(1+\|Y\|^T_{\infty})^2}\mathbb{E}\left[\int_t^{T}|Z_s|^2ds\Bigg|\mathcal{F}_t\right]\leq\mathbb{E}\left[\int_t^{T}\frac{1}{2(1+Y_s)^2}|Z_s|^2ds\Bigg|\mathcal{F}_t\right]\\
&\leq \|Y\|^T_{\infty}\leq C<+\infty,
\end{align*}
thus $Z\in\text{BMO}(\mathbb{P}).$
\end{proof}
\begin{proof}[Proof of Corollary \ref{cor:Linff}]
The BSDE
$$Y^h_t=X+\int_t^T \left(\alpha_sY^h_s+\beta_sY^h_s|\ln(Y^h_s)|+\gamma_s\cdot Z^h_s+\delta|Z^h_s|^2/Y^h_s \right)ds-\int_t^TZ^h_sdW_s$$
can be rewritten as
\begin{equation}
Y^h_t=X+\int_t^T \left(\alpha_sY^h_s+\beta_sY^h_s|\ln(Y^h_s)|+\delta|Z^h_s|^2/Y^h_s \right) ds-\int_t^TZ^h_sdW^{\gamma}_s,
\label{eq:proof1}
\end{equation}
where Girsanov's theorem is employed. 
Then, by Proposition~\ref{prop:Linf}, there exists a (unique) solution $(Y^h,Z^h)$ to Equation~\eqref{eq:proof1} such that $|Y^h_t|\leq K$  $d\mathbb{Q}^{\gamma}\times dt\text{-a.s.}$ and $Z^h\in\text{BMO}(\mathbb{Q}^{\gamma})$. 
Here, $\mathbb{Q}^{\gamma}$ is the probability measure with density $\mathcal{E}^{\gamma}_T$. 
As is well known, $\mathbb{Q}^{\gamma}$-boundedness implies $\mathbb{P}$-boundedness. 
Furthermore, any $\text{BMO}(\mathbb{Q}^{\gamma})$ martingale is also a $\text{BMO}(\mathbb{P})$ martingale (see for instance \cite{BEK05}). 
Hence, Equation~\eqref{eq:proof1} admits a unique solution $(Y^h,Z^h)\in\mathcal{H}^{\infty}_T\times\text{BMO}(\mathbb{P}).$ 
The existence of a solution $(Y,Z)$ with the required regularity for the general BSDE with parameters $(X,g)$ follows similarly as in the proofs of Theorem~\ref{th:Egen} and Proposition~\ref{prop:Linf}. 
\end{proof}

\begin{proof}[Proof of Lemma~\ref{lem:biharine}]
We want to apply Theorem~1 in \cite{HKK22} to the function 
$$\bar\psi(x)=\begin{cases}
2\ln(2)\ &\text{ if } 0\leq x\leq2, \\
x\ln(x) &\text{ if } x>2.
\end{cases}$$
It is clear that $\bar\psi$ is convex, positive and non-decreasing on $\R_+$. 
Furthermore, we have that $x|\ln(x)|\leq\bar\psi(x)$ for any $x\in\R_+$ and $\psi(x)=\int_2^x\frac{1}{\bar\psi(r)}dr,$ for any $x>0$. 
The integrability conditions on $X,\beta,u$ guarantee that $\mathbb{E}[|X+\int_0^T\beta_s\bar\psi(u_s)ds|^{p'e^B}]<+\infty,$ for some $p'>1$. 
Indeed, ${\bar\psi(x)}\leq C_{\varepsilon}+ x^{\varepsilon}$ for any $\varepsilon>1$, with $C_{\varepsilon}>0$. 
Thus, it holds that
 $$\mathbb{E}\left[X+\int_0^T\beta_s\bar\psi(u_s)ds\right]\leq C_{\varepsilon}+\mathbb{E}\left[X+\int_0^T\beta_su_s^{\varepsilon}ds\right]\leq C_{\varepsilon}+\mathbb{E}\left[X+\|\beta\|^T_{\infty}\esssup_{t\in[0,T]}u_s^{\varepsilon}\right].$$ 
Choosing $\varepsilon>1$ such that $p'e^B=\varepsilon pe^B=p(e^B+1)$ (hence, $\varepsilon=1+e^{-B}$) with $p'=\varepsilon p>1$, the integrability follows.
Thus, we only need to check that 
$$\psi^{-1}\left(\psi\left(\mathbb{E}\left[X+\int_t^T\beta_s\bar\psi(u_s)ds\right]\right)+\int_0^t\beta_sds\right)$$ is uniformly integrable in time.
Setting $X_t:=\mathbb{E}\left[X+\int_t^T\beta_s\bar\psi(u_s)ds\Big|\mathcal{F}_t\right]$ it holds by direct inspection that
\begin{align*}
&\psi^{-1}\left(\psi(X_t)+\int_0^t\beta_sds\right)=\exp\left(\exp(\ln(\ln(X_t))+\int_0^t\beta_sds+\ln(\ln(2)))\right) \\
&\leq\exp(e^B\ln(X_t)))=(X_t)^{e^B}=\mathbb{E}\left[X+\int_t^T\beta_s\bar\psi(u_s)ds\bigg|\mathcal{F}_t\right]^{e^B} \\
&\leq\mathbb{E}\left[X+\int_0^T\beta_s\bar\psi(u_s)ds\bigg|\mathcal{F}_t\right]^{e^{B}}.
\end{align*}
The integrability conditions on $X,\beta,u$ and the properties of conditional expectations ensure that 
$$\mathbb{E}\left[\Bigg|\mathbb{E}\left[X+\int_0^T\beta_s\bar\psi(u_s)ds\bigg|\mathcal{F}_t\right]\Bigg|^{p'e^B}\right]\leq\mathbb{E}\left[\Bigg|X+\int_0^T\beta_s\bar\psi(u_s)ds\Bigg|^{p'e^B}\right]<+\infty,$$ thus the process $\psi^{-1}(\psi(X_t)+\int_0^t\beta_sds)$ is bounded in $L^{p'}$, $p'>1$. 
Hence, it is uniformly integrable and the thesis follows.
\end{proof}

\begin{proof}[Proof of Corollary \ref{cor:stability}]
As is clear from Theorem~\ref{th:stability}, we only need to verify the convergence $\int_0^Tg^n(t,Y_t,Z_t)dt\to \int_0^Tg(t,Y_t,Z_t)dt$ in $L^p$ for any $p\geq 1$ to establish the thesis. 
For any fixed $p\geq 1$, it holds that
\begin{align*}
&\mathbb{E}\left[\left|\int_0^Tg^n(t,Y_t,Z_t)dt\right|^p\right]\leq\mathbb{E}\left[\left(\int_0^T|g^n(t,Y_t,Z_t)|dt\right)^p\right] \\
& \leq\mathbb{E}\left[\left(\int_0^T\left(\alpha_t Y_t+\beta_t Y_t|\ln(Y_t)|+\gamma_t|Z_t|+\delta \frac{|Z_t|^2}{Y_t}\right)dt\right)^p\right].
\end{align*}
Then, by the assumptions $X\in L^p$ for any $p\geq 1$, $X\geq \varepsilon$ and applying Proposition~\ref{prop:Zreg}, we obtain uniform integrability of $\int_0^Tg^n(t,Y_{t},Z_{t})dt$ for any $p\geq 1$. 
In addition, $g^n(t,V_t,Z_t)\xrightarrow{n\to\infty}g(t,V_t,Z_t)$ $d\mathbb{P}\times dt$-a.s., thus Vitali's theorem leads to the thesis.
\end{proof}

\begin{proof}[Proof of Proposition~\ref{prop:2DU}]
Consider two solutions $(Y,Z)$ and $(Y',Z')$ to Equation~\eqref{eq:twodrivers} such that $0< Y,Y'\leq Y^{h_1,h_2}$. 
Then the corresponding one-driver equation with parameters $(X,g)$ with $g(t,y,v):=g_1(t,y,g_2^{-1}(t,y,v))$ verifies all the assumptions of Theorem~\ref{th:EUg1}, hence there exists a unique solution $(Y^{g},Z^{g})$ to the BSDE with parameters $(X,g)$. 
Clearly, this implies the uniqueness of the first component of the solution to Equation~\eqref{eq:twodrivers}. 
We only need to verify the uniqueness for the $z$-component. 
Proceeding exactly as in the proof of Theorem~\ref{th:EUg1}, we can prove that
$$\int_0^T|g_2(t,Y_t,Z_t)-g_2(t,Y_t,Z'_t)|^2=0 \ \text{ a.s.,}$$ which entails
$$g_2(t,Y_t,Z_t)=g_2(t,Y_t,Z_t'), \ \  d\mathbb{P}\times dt\text{-a.s.}$$ 
Injectivity of $g_2$ w.r.t.\ $z$ gives $Z_t=Z_t'$ $d\mathbb{P}\times dt$-a.s.
\end{proof}

\begin{proof}[Proof of Corollary~\ref{cor:euGBSDE}]
Let us observe that Equation~\eqref{eq:GBSDE} can be rewritten as a two-driver BSDE, by defining $g_1(t,y,z):=y\tilde f(t,y,z)$ and $g_2(t,y,z):=yz$. 
Hence, the part of the statement regarding existence is a straightforward consequence of Proposition~\ref{prop:twodriverE}. 
For the uniqueness, we have that $g_1(t,y,g_2^{-1}(t,y,v))=y\tilde f(t,y,v/y)=y\tilde f_1(t,y)+y \tilde f_2(t,v/y).$ 
Thus, $y\mapsto y\tilde f_1(\fatdot{}, y)$ is convex by assumption and $(y,v)\mapsto y\tilde f_2(v/y)$ is jointly convex in $(y,v)$ being the perspective function of the convex function $\tilde f_2$. 
Hence, $g_1(t,y,g_2^{-1}(t,y,v))$ is jointly convex in $(y,v)$ and uniqueness follows from Proposition~\ref{prop:2DU}. 
The last thesis is obvious, by taking $\tilde f_1\equiv0.$
\end{proof}

\begin{proof}[Proof of Proposition~\ref{prop:return_gbsde}]
Existence and uniqueness follow from Corollary~\ref{cor:euGBSDE}.
Monotonicity and the Lebesgue property follow from Propositions~\ref{prop:2Dcomparisongen}
and \ref{prop:2Dstability} (and Remark~\ref{rem:stabless}), respectively.

For GA-convexity, we first work on $\mathcal{L}^{\infty}(\mathcal{F}_T)$.
Setting $f(t,y,z)=\tilde f(t,e^y,z)+\frac{1}{2}|z|^2$, the GA-convexity assumption on
$\tilde f$ implies convexity of $f$ in $(y,z)$, which by the BSDE comparison theorem yields
convexity of the corresponding monetary risk measure $\rho_t(\cdot)=\tilde\rho_t(e^{\cdot})$.
Multiplicative convexity of $\tilde\rho_t$ on $\mathcal{L}^{\infty}$ then follows from
Proposition~\ref{prop:1to1P}~(e). The extension to $L_+^{p(2\delta+1)(e^B+1)}$ uses a
denseness argument: setting $X^n:=n\wedge X\vee\frac{1}{n}\in\mathcal{L}^{\infty}$, one has
$X^n\to X$ in $L^{p(2\delta+1)(e^B+1)}$ by dominated convergence; the stability result
(Proposition~\ref{prop:2Dstability}) then allows passage to the limit a.s.\ along subsequences.

Positive homogeneity follows because $f$ independent of $y$ gives cash-additivity of
$(\rho_t)$, hence positive homogeneity of $(\tilde\rho_t)$ by Proposition~\ref{prop:1to1P}~(b).
If $\tilde f$ is increasing in $y$, the same monotonicity carries over to $f$, yielding
cash-superadditivity and hence star-shapedness via Proposition~\ref{prop:1to1P}~(c).
Both properties extend to $L_+^{p(2\delta+1)(e^B+1)}$ by the same denseness argument.

Time-consistency on $\mathcal{L}^{\infty}$ follows from the flow property of BSDEs
(\cite{EPQ97}) and Proposition~\ref{prop:1to1P}~(g), and extends by denseness.
Normalization is verified by direct inspection: $(\tilde\rho,\tilde Z)\equiv(1,0)$ solves
\eqref{eq:genGBSDEs} with $X\equiv 1$ whenever $\tilde f(\cdot,1,0)\equiv0$, and uniqueness
concludes.
\end{proof}


{\footnotesize
\begin{thebibliography}{999}
\setlength{\itemsep}{1pt} 
\setlength{\parskip}{1pt} 

\bibitem{AS08}
\textsc{Aumann, R.~J. and R. Serrano} (2008). 
An economic index of riskiness. 
\textit{Journal of Political Economy} 116(5), 810-836.

\bibitem{B20} \textsc{Bahlali, K.} (2020). 
A domination method for solving unbounded quadratic BSDEs.
\textit{Graduate Journal of Mathematics} 5, 20-36.

\bibitem{BEO17} \textsc{Bahlali, K., M. Eddahbi and Y. Ouknine} (2017). 
Quadratic {BSDE} with {$\mathbb{L}^2$}-terminal data: {K}rylov's estimate, {I}t\^{o}-{K}rylov's formula and existence results.
\textit{The Annals of Probability} 45(4), 2377-2397.

\bibitem{BEH15} \textsc{Bahlali, K., E. Essaky and M. Hassani} (2015). 
Existence and uniqueness of multidimensional BSDEs and of systems of degenerate PDEs with superlinear growth generator.
\textit{SIAM Journal on Mathematical Analysis} 47(6), 4251-4288.

\bibitem{BKK17} \textsc{Bahlali, K., O. Kebiri, N. Khelfallah and H. Moussaoui} (2017). 
One dimensional BSDEs with logarithmic growth application to PDEs.
\textit{Stochastics} 89, 1061-1081.

\bibitem{BT19} \textsc{Bahlali, K. and L. Tangpi} (2021). 
BSDEs driven by $|z|^2/y$ and applications to PDEs and decision theory.
\textit{Preprint on arXiv: https://arxiv.org/abs/1810.05664v3}.

\bibitem{BEK05} \textsc{Barrieu, P. and N. El Karoui} (2005). 
Pricing, hedging and optimally designing derivatives via minimization of risk measures. 
In: Carmona, R. (ed.). \textit{Indifference Pricing: Theory and Applications}, 77-144,
Princeton University Press, Princeton.

\bibitem{BZ16} \textsc{Bauer, D. and G. Zanjani} (2016).
The marginal cost of risk, risk measures, and capital allocation.
\textit{Management Science} 62(5), 1431-1457.

\bibitem{BLR18} \textsc{Bellini, F., R.~J.~A. Laeven and E. Rosazza Gianin} (2018).
Robust return risk measures. 
\textit{Mathematics and Financial Economics} 12, 5-32.

\bibitem{BLR21} \textsc{Bellini F., R.~J.~A. Laeven and E. Rosazza Gianin} (2021).
Dynamic robust Orlicz premia and Haezendonk-Goovaerts risk measures. 
\textit{European Journal of Operational Research} 291, 438-446.

\bibitem{BN08} \textsc{Bion-Nadal, J.} (2008). 
Dynamic risk measures: Time consistency and risk measures from BMO martingales. 
\textit{Finance and Stochastics} 12, 219-244.

\bibitem{BH06} \textsc{Briand, P. and Y. Hu} (2006). 
BSDE with quadratic growth and unbounded terminal value.
\textit{Probability Theory and Related Fields} 136, 604-618.

\bibitem{BH08} \textsc{Briand, P. and Y. Hu} (2008). 
Quadratic BSDEs with convex generators and unbounded terminal conditions.
\textit{Probability Theory and Related Fields} 141, 543-567.

\bibitem{CCMTW22} \textsc{Castagnoli, E.,  G. Cattelan, F. Maccheroni, C. Tebaldi and R. Wang} (2022). 
Star-shaped risk measures. 
\textit{Operations Research} 70(5), 2637-2654.

\bibitem{CF10} \textsc{Chateauneuf, A. and J.~H. Faro} (2010). 
Ambiguity through confidence functions. 
\textit{Journal of Mathematical Economics} 45, 535-558.

\bibitem {CE02} \textsc{Chen, Z. and L.~G. Epstein} (2002).
Ambiguity, risk, and asset returns in continuous time.
\emph{Econometrica} 70, 1403-1443.

\bibitem{CHMP} \textsc{Coquet, F., Y. Hu, J. M\'emin and S. Peng} (2002). 
Filtration-consistent nonlinear expectations and related $g$-expectations. 
\textit{Probability Theory and Related Fields} 123(1), 1-27.

\bibitem{D12} \textsc{Delbaen, F.} (2012). 
\textit{Monetary Utility Functions}.
Osaka University Press, Osaka.

\bibitem{DHR11} \textsc{Delbaen, F., Y. Hu and A. Richou} (2011).
On the uniqueness of solutions to quadratic BSDEs with convex generators and unbounded terminal conditions.
\textit{Annales de l’Institut Henri Poincaré -- Probabilités et Statistiques} 47(2), 559-574.

\bibitem{DPR10} \textsc{Delbaen, F., S. Peng and E. Rosazza Gianin} (2010).
Representation of the penalty term of dynamic concave utilities. 
\textit{Finance and Stochastics} 14, 449-472.

\bibitem{DR09} \textsc{Detemple, J. and M. Rindisbacher} (2009). 
Dynamic asset allocation: Portfolio decomposition formula and applications. 
\textit{The Review of Financial Studies} 23(1), 25-100.

\bibitem{DHK13} \textsc{Drapeau, S., G. Heyne and M. Kupper} (2013). 
Minimal supersolutions of convex BSDEs. 
\textit{The Annals of Probability} 41, 3973-4001.

\bibitem{DK12} \textsc{Drapeau, S. and M. Kupper} (2013). 
Risk preferences and their robust representation. 
\textit{Mathematics of Operations Research} 38, 28-62.

\bibitem{DKRT14} \textsc{Drapeau, S., M. Kupper, E. Rosazza Gianin and L. Tangpi} (2016). 
Dual representation of minimal supersolutions of convex BSDEs. 
\textit{Annales de l'Institut Henri Poincaré -- Probabilités et Statistiques} 52, 868-887.

\bibitem {DE92} \textsc{Duffie, D. and L.~G. Epstein} (1992).
Stochastic differential utility.
\textit{Econometrica} 60, 353-394.

\bibitem{EPQ97} \textsc{El Karoui, N., S. Peng and M.~C. Quenez} (1997).
Backward stochastic differential equations in finance. 
\textit{Mathematical Finance} 7, 1-71.

\bibitem{ELKR09} \textsc{El Karoui, N. and C. Ravanelli} (2009). 
Cash subadditive risk measures and interest rate ambiguity. 
\textit{Mathematical Finance} 19, 561-590.

\bibitem{EK86} \textsc{Ethier, S.~N. and T.~G. Kurtz} (1986). 
\textit{Markov Process: Characterization and Convergence}.
J. Wiley, New York.

\bibitem{FHT23} \textsc{Fan, S., Y. Hu and S. Tang} (2023). 
Existence, uniqueness and comparison theorem on unbounded solutions of scalar super-linear BSDEs.
\textit{Stochastic Processes and Their Applications} 157, 335-375.

\bibitem{FS11} \textsc{F\"ollmer, H. and A. Schied} (2011). 
\textit{Stochastic Finance}.
3rd ed., De Gruyter, Berlin.

\bibitem{FH09} \textsc{Foster, D.~P. and S. Hart} (2009). 
An operational measure of riskiness. 
\textit{Journal of Political Economy} 117(5), 785-814.

\bibitem{HIM05} \textsc{Hu, Y., P. Imkeller and M. M\"uller} (2005). 
Utility maximization in incomplete markets.
\textit{The Annals of Applied Probability} 15(3), 1691-1712.

\bibitem{HKK22} \textsc{Hun, O., M.-C. Kim and K.-G. Kim} (2022). 
Stochastic Bihari inequality and applications to BSDE.
\textit{Journal of Mathematical Analysis and Applications} 513(1), 126204.

\bibitem{J08} \textsc{Jiang, L.} (2008). 
Convexity, translation invariance and subadditivity for $g$-expectations and related risk measures. 
\textit{The Annals of Applied Probability} 18(1), 245-258.

\bibitem{K94} \textsc{Kazamaki, N.} (1994).
\textit{Continuous Exponential Martingales and BMO}. Lecture Notes in Mathematics 1579, Springer, Berlin.

\bibitem{K00} \textsc{Kobylanski, M.} (2000).
Backward stochastic differential equations and partial differential equations with quadratic growth.
\textit{The Annals of Probability} 28(2), 558-602.

\bibitem{KLLSS18} \textsc{Kr\"atschmer, V., M. Ladkau, R.~J.~A. Laeven, J.~G.~M. Schoenmakers and M.~A. Stadje} (2018). 
Optimal stopping under uncertainty in drift and jump intensity. 
\textit{Mathematics of Operations Research} 43, 1177-1209.

\bibitem{LR22} \textsc{Laeven, R.~J.~A. and E. Rosazza Gianin} (2022). 
Quasi-logconvex measures of risk. 
\textit{Preprint on arXiv:2208.07694v1}.

\bibitem{LRZ24} \textsc{Laeven, R.~J.~A., E. Rosazza Gianin and M. Zullino} (2024). 
Law-invariant return and star-shaped risk measures.
\textit{Insurance: Mathematics and Economics} 117, 140-153.

\bibitem{LRZ23} \textsc{Laeven, R.~J.~A., E. Rosazza Gianin and M. Zullino} (2026). 
Star-shaped and dynamic return risk measures via BSDEs.
\textit{Finance and Stochastics} 30, 903-950.

\bibitem{LSSS24}
\textsc{Laeven, R.~J.~A., J.~G.~M. Schoenmakers, N.~F.~F. Schweizer and M.~A. Stadje} (2025). 
Robust multiple stopping --- A duality approach. 
\textit{Mathematics of Operations Research}, 50, 1250-1276.

\bibitem{LS13} \textsc{Laeven, R.~J.~A. and M.~A. Stadje} (2013). 
Entropy coherent and entropy convex measures of risk. 
\textit{Mathematics of Operations Research} 38, 265-293.

\bibitem{LS14} \textsc{Laeven, R.~J.~A. and M.~A. Stadje} (2014). 
Robust portfolio choice and indifference valuation.  
\textit{Mathematics of Operations Research} 39, 1109-1141.

\bibitem{LQ03} \textsc{Lazrak, A. and M.~C. Quenez} (2003). 
A generalized stochastic differential utility. 
\textit{Mathematics of Operations Research} 28, 154-180.

\bibitem{LMO22} \textsc{Likibi Pellat, R.,  
O. Menoukeu Pamen and Y. Ouknine} (2022).
A class of quadratic forward-backward stochastic differential equations.
\textit{Journal of Mathematical Analysis and Applications} 514, 126100.

\bibitem{MS05} \textsc{Mania, M. and M. Schweizer} (2005). 
Dynamic exponential utility indifference valuation. 
\textit{The Annals of Applied Probability} 15(3), 2113-2143.

\bibitem{PDS23} \textsc{Papapantoleon, A., D. Possamaï and A. Saplaouras} (2023).
Stability of backward stochastic differential equations: The general Lipschitz case.
\textit{Electronic Journal of Probability} 28, 1–56.

\bibitem{PP90} \textsc{Pardoux, E. and S. Peng} (1990). 
Adapted solution of a backward stochastic differential equation. 
\textit{Systems and Control Letters} 14, 55-61.

\bibitem{P97} \textsc{Peng, S.} (1997). 
Backward SDE and related $g$-expectation.
In: El Karoui, N. and L. Mazliak (eds.). 
\textit{Backward Stochastic Differential Equations}, 141-159,
Longman, Harlow.

\bibitem{P05} \textsc{Peng, S.} (2005). 
Dynamically consistent nonlinear evaluations and expectations. 
\textit{Preprint on arXiv:math.PR/0501415}.

\bibitem{P64} \textsc{Pratt, J.~W.} (1964).
Risk aversion in the small and in the large.
\textit{Econometrica} 32, 122-136.

\bibitem{Pr05} \textsc{Protter, P.} (2005). 
\textit{Stochastic Integration and Differential Equations}.
Springer Science \& Business Media.

\bibitem{RG06} \textsc{Rosazza Gianin, E.} (2006). 
Risk measures via $g$-expectations. 
\textit{Insurance: Mathematics and Economics} 39, 19-34.

\bibitem{W08} \textsc{Wakker, P.~P.} (2008).
Explaining the characteristics of the power (CRRA) utility family. 
\textit{Health Economics} 17, 1329-1344.

\bibitem{WF18} \textsc{Wang, X. and S. Fan} (2018).
A class of stochastic Gronwall's inequality and its application.
\textit{Journal of Inequalities and Applications} 336, 1-10.

\bibitem{Z17} \textsc{Zhang, J.} (2017). 
\textit{Backward Stochastic Differential Equations}.
Springer, New York, NY.
\end{thebibliography}
}
\end{document}